\documentclass[10pt,a4paper,twoside]{article}
\usepackage{layout}
\usepackage{exscale, amsmath, amsfonts, amssymb, amsthm, amscd}
\usepackage{verbatim}
\newcommand{\MMM}{}

\newcommand{\ov}[1]{{\overline #1}}

\newcommand{\convinfty}[1]
  {\stackrel{#1\to\infty}{-\hspace{-2mm}
  -\hspace{-2mm}-\hspace{-2mm}\longrightarrow}}

\newcommand{\Area}{{\rm area\,}}
\newcommand{\Inter}{{\rm int\,}}
\newcommand{\dist}{{\rm dist}}

\newcommand{\N}{\mathbb{N}}

\newcommand{\R}{\mathbb{R}}
\newcommand{\Z}{\mathbb{Z}}
\newcommand{\T}{\mathbb{T}}

\newcommand{\w}{\omega}
\newcommand{\W}{\Omega}

\newcommand{\A}{{\mathcal A}}
\newcommand{\B}{{\mathcal B}}
\newcommand{\C}{{\mathcal C}}
\newcommand{\F}{{\mathcal F}}

\newcommand{\I}{{\mathcal I}}
\newcommand{\J}{{\mathcal J}}
\renewcommand{\L}{{\mathcal L}}
\newcommand{\M}{{\mathcal M}}
\renewcommand{\P}{{\mathcal P}}
\newcommand{\Tail}{{\mathcal T}}

\renewcommand{\t}{\tau}

\newcommand{\ErgSum}[2]{\sum_{i=0}^{n-1}#1\circ #2^i}

\newcommand{\Lc}{{\widehat{L}}}
\newcommand{\Lf}{{\widehat{L}}^\sharp}

\newcommand{\len}{{\rm length}\,}
\newcommand{\area}{{\rm area}\,}
\newcommand{\inter}{{\rm int}\,}

\newcommand{\nd}{\noindent}

\newtheorem{thm}{{\bf  Theorem}}[section]
\newtheorem{lemma}[thm]{{\bf  Lemma}}
\newtheorem{cor}[thm]{{\bf  Corollary}}

\newtheorem{definition}[thm]{{\bf  Definition}}
\newtheorem{rem}[thm]{{\bf Remark}}

\newcommand{\pf}{\noindent {\bf Proof. }}
\newcommand{\pfofthm}{\noindent {\bf Proof of the theorem. }}

\newcommand{\nn}{\nonumber}

\def\blfootnote{\xdef\@thefnmark{}\@footnotetext}
\long\def\symbolfootnote[#1]#2{\begingroup%
\def\thefootnote{\fnsymbol{footnote}}\footnote[#1]{#2}\endgroup}

\renewcommand{\thefootnote}{\symbol{footnote}}

\begin{document}
%>>>>>>>>>>>>>>>>>>>>>>>>>>>>>>>>>>>>>>>>>>>>>>>>>>>>>>>>>>>>>>>>>>>
%\end{document}
%TTTTTTTTTTTTTTTTTTTTTTTTTTTTTTTTTTTTTTTTTTTTTTTTTTTTTTTTT
% Title
%TTTTTTTTTTTTTTTTTTTTTTTTTTTTTTTTTTTTTTTTTTTTTTTTTTTTTTTTT

\thispagestyle{empty}

\vspace{1mm}

\nd\rule{6in}{1pt}

\smallskip

\nd {\large J.\,Brettschneider}

%\symbolfootnote[0]{
%J.\,Brettschneider: Dept.~of Statistics, University of Warwick, Coventry, CV4 7AL, UK;\\
%\phantom{I23} 
%Dept.~of Community Health $\&$ Epidemiology, Cancer Research Institute Div.~of Cancer Care $\&$ Epidemiology,\\
%\phantom{I23}  Queen's University, Kingston, Ontario, K7L 3N6, Canada .\\
%\phantom{12.3}{\tt julia.brettschneider@warwick.ac.uk}}

\smallskip

\nd {\large\bf Shannon-McMillan theorems for discrete random fields along curves} \\
{\large\bf and lower bounds for surface-order large deviations}\\

%\vspace{1mm}
\thispagestyle{plain}

%ABABABABABABABABABABABABABABABABABABABABABABABBABABABAB
%             Abstract
%ABABABABABABABABABABABABABABABABABABABABABABABBABABABAB

\nd {\bf Abstract:}
The notion of a surface-order specific entropy $h_c(P)$ 
of a two-dimensional discrete random field $P$ along a curve $c$ 
is introduced as the limit of rescaled entropies along lattice 
approximations of the blowups of $c.$ 
Existence is shown by proving a corresponding Shannon-McMillan theorem. 
We obtain a representation of $h_c(P)$ as a mixture of 
specific entropies along the tangent lines of $c.$  
As an application, the specific entropy along curves is used to
refine F{\"o}llmer and Ort's lower bound for the large deviations
of the empirical field of an attractive Gibbs measure from its ergodic
behaviour in the phase-transition regime.

\smallskip

\nd\rule{6in}{1pt}

%SSSSSSSSSSSSSSSSSSSSSSSSSSSSSSSSSSSSSSSSSSSSSSSSSSSSSSSSSSSSSSSS
\section{Introduction}\label{int}
%\markboth{Skew products}{Introduction}
%SSSSSSSSSSSSSSSSSSSSSSSSSSSSSSSSSSSSSSSSSSSSSSSSSSSSSSSSSSSSSSSS

The entropy of a stationary random field $P$ indexed by a 
lattice is usually defined as a limit of entropies on an
increasing sequence of boxes, rescaled by the volume of the boxes.
Shannon-McMillan theorems describe this convergence on a deeper
level, as $\L^1$-convergence of rescaled information quantities. 
In the context of large deviations for Gibbs measures the volume-order 
entropies may not provide enough information when a phase transition 
occurs.  Instead, F{\"o}llmer and Ort \cite{FoOrt88} 
introduced the concept of surface-order entropy on boxes,
proved a corresponding version of the Shannon-McMillan theorem 
and used it to estimate large-deviation probabilities.
The construction of the Wulff shape by Dobrushin,
Kotecky, and  Shlosman \cite{DKS92} suggests that such estimates
can be improved if boxes are replaced by more general shapes.

We investigate the problem of constructing entropies on general surfaces, 
proving appropriate versions of the Shannon-McMillan theorem,
and using these constructions to refine the large deviation lower bound. 
We carry out this program in the two-dimensional 
case, where surfaces reduce to contour curves. 
The existence of surface-order specific entropy does not simply 
follow from a subadditivity argument. Instead, we prove a corresponding
Shannon-McMillan theorem, that is, $L^1(P)$-convergence of rescaled
information quantities along lattice approximations of the successive blowups 
of the curve $c.$
This is accomplished in three steps.  The first main result is a
Shannon-McMillan theorem along lines;  cf.~Theorem \ref{SMQ} for 
rational slopes and Theorem \ref{SMnotQ} for irrational slopes.
The proof relies on {\it uniform} convergence in ergodic theorems  
for a suitable skew product transformation. 
The second step is to extend the results to polygons under the
assumption of a strong $0$-$1$ law on $P.$  Finally, we obtain
an explicit formula for the specific entropy $h_{c}(P)$ as a mixture of
the conditional entropy of the random field restricted to the origin, given 
suitably defined {\it past}-$\sigma$-algebras along tangent lines of the curve
(Theorem \ref{ShMcMCurve}).
Under certain conditions, this construction can be extended to relative 
entropies of one random field with respect to another. This will be the key 
to our proof of a refined lower bound for large deviations in the
phase transition regime (Theorem \ref{LowerBd}).

The role of Shannon-McMillan theorems in the refined analysis
of large deviations provided the original motivation for this work.
It seems, however, that the study of entropies
along surfaces may hold independent interest.
This paper lays some of the groundwork for such
a general theory of specific entropies along shapes.  
The key to this theory is a careful combination of probabilistic 
arguments and non-standard discrete geometrical constructions. 
We now explain our approach and our results in more detail.

{\bf Shannon-McMillan theorems and entropy.}
Consider a random sequence $\w$
%:\{1,2,\dots\}\longrightarrow\Upsilon$
of letters from a finite alphabet $\Upsilon,$ modelled by a
probability measure $P$ on the space
$\W:=\Upsilon^{\{1,2,\dots\}}.$ For any finite $n,$ the {\it
information} provided by the first $n$ letters can be described by
the function $ -\log P[\,\w_{\{1,\dots,n\}}], $ where
$\w_{\{1,\dots,n\}}$ denotes the restriction of $\w$ to
$\{1,\dots,n\},$ and $P\big[(y_1,\dots,y_n)\big]$ is the
probability that the pattern $(y_1,\dots,y_n)$ appears in the
first $n$ trials. In the classical case, when the letters are
independent and identically distributed according to a measure
$\mu,$ the classical Shannon-McMillan theorem states that the
rescaled information functions $ -n^{-1}\log P[\,\w_{\{1,\dots,n\}}], $ 
converge in $\L^1(P)$ to the {\it entropy} 
$ H(\mu):=-\sum_{y\in\Upsilon}\mu(y)\log\mu(y) $ of the
measure $\mu.$ The theorem can be extended to a general ergodic
sequence.  In the bilateral case, when $P$ is an ergodic measure
on $\Upsilon^{\Z},$ the limiting quantity takes the form 
$h(P):=E\big[H\big(P_0[\,\cdot\,\vert\,\P](\w)\big)\big], $ where
$P_0[\,\cdot\,\vert\,\P]$ is the distribution of the random field in $(0,0)$
conditioned on the``past'' $\sigma$-algebra  $\P$ 
%$F_{\{-1,-2,\dots\}}$
generated by the projection of $\w$ to the set $\{-1,-2,\dots\}.$

These constructions can be extended to a spatial setting when
the random field is given by a stationary probability measure $P$
on a configuration space $\Upsilon^{\Z^{d}}.$
Thouvenot \cite{Tho72} and F{\"o}llmer \cite{Foe73} proved
spatial extensions of Shannon and McMillan's result.
The {\it specific entropy} $h(P)$ is introduced as the limit of 
$\vert V_n\vert^{-1}H_{V_n}(P),$
where $V_n$ is the set of all lattice sites in the box $[-n,n]^d,$ and 
$H_{V_n}(P)$ is the entropy of the measure $P$ restricted to $V_n.$
The existence of the limit follows from the subadditivity of
$H_{V}$ with respect to $V.$
What is more, the corresponding Shannon-McMillan theorem shows
$\L^1(P)$-convergence of the rescaled information functions
$-\vert V_n\vert^{-1}\log P[\w_{V_n}].$
If $P$ is ergodic we obtain the formula
$h(P)=E\big[H\big(P_0[\,\cdot\,\vert\,\P^d](\w)\big)\big],$
where $\P^d$ is a $\sigma$-algebra representing a
spatial version of the ``past''. More precisely, $\P^d$
is generated by the projections of $\w$ to the sites
preceding the origin in the lexicographical ordering of $\Z^d.$

{\bf Surface entropy.} 
Our goal is to derive refined versions
of the Shannon-McMillan theorem, where the information functions
are observed along surfaces. This was carried out in \cite{FoOrt88} 
for the surfaces of boxes parallel to the axes.
In this paper, we focus on the two-dimesional case. We develop a
construction of surface entropy where rectangles are replaced by
general curves. More precisely, following a suggestion by 
Hans F{\"o}llmer, we introduce the specific entropy along 
sets generated by lattice approximations to lines, and then
extend this to polygons and piecewise smooth curves.  

Most of the work here goes into our first result,
a Shannon-McMillan theorem for the {\it specific entropy} $h_{\lambda}(P)$
of a stationary random field $P$ {\it along a line}
% is the foundation of this work 
(Theorem \ref{SMQ} if the slope is rational and Theorem \ref{SMnotQ} 
if the slope is irrational). 
We prove the $\L^1(P)$-convergence of the rescaled information
functions along increasing segments of the  line's {\it lattice approximation}
$\{(z,[\lambda z+a])\,\vert\,z\in\Z\},$
where $[x]$ denotes the integer part of $x.$
If $P$ fulfills a $0$-$1$ law on the tail field, we obtain the formula
\begin{equation*}%\label{EntrLine}
h_{\lambda}(P)=\int_0^1 E \big[H\big(P_0[\,\cdot\,
\vert\,\P_{\lambda,t}](\w)\big)\big]\,dt,
\end{equation*}
where $\P_{\lambda,t}$ is the $\sigma$-algebra
generated by those approximating sites
which precede $0$ in the lexicographical ordering of $\Z^2.$  
If $\lambda$ is rational, the it suffices for $P$ to be ergodic.
Furthermore, the specific entropy along the line can be written as
$$
h_{\lambda}(P)=
\frac{1}{q}\sum_{\nu=0}^{q-1}E\big[H\big(P_{0}\big[\cdot\,\vert
\P_{{p/q,\,\nu p/q}}\big](\w)\big)\big], 
$$
where $p/q$ is the unique representation of $\lambda$ by
integers $p\in\N$ and $q\in\Z$ having no common divisor. The past
$\sigma$-algebras $\P_{p/q,\,\nu p/q}$ correspond to
the $q$ different possibilities to start the $q$-periodic pattern
of the lattice approximation.

The idea to investigate such an specific entropy along lines 
has two precursors. The first is a volume-order {\it directional entropy}, 
which Milnor (\cite{Mil86} and \cite{Mil88}) introduced in the 
context of cellular automata.  The second is the specific entropy
along hyperplanes perpendicular to an axis, which was defined 
by \cite{FoOrt88} as a step toward their surface entropy along boxes.
It may be noted that a distinction between rational and 
irrational slopes was also made by Sinai in his work \cite{Sin85} on
Milnor's directional entropy for cellular automata.
This construction was further developed by Park (\cite{Par94}, 
\cite{Par95}, and \cite{Par96}) and Sinai \cite{Sin94}.
The original problem of continuity with respect to the
direction was finally solved in Park \cite{Par99}.

Our next result is a Shannon-McMillan theorem along polygons
(Theorem \ref{ShMcMpolygon}). In particular, we obtain a
representation of the {\it specific entropy} of $P$
{\it along a polygon} as a mixture of entropies
along lines corresponding to its edges. 
The extension to polygons requires a strong form
of the $0$-$1$ law on the tail field, which was introduced
in \cite{FoOrt88}. It says that, for any subset $J$ of
$\Z^2,$ the $\sigma$-algebra generated by the sites in $J$
does not increase if we add information about the tail
behaviour outside of $J;$ cf.~Definition \ref{stzeroone}.
Finally, we come to the main result in this section.
Theorem \ref{ShMcMCurve} says that the {\it specific entropy along a curve}
$c:[0,T]\longrightarrow\R^2$ is a mixture of entropies along the tangent lines:
\begin{equation}\label{MainResult}
h_c(P)=\int_0^T h_{c'(t)}(P)\,dt
\end{equation}
Here, $h_{c'(t)}(P)$ denotes the specific entropy along a line having the same
slope as the tangent of $c$ in $t;$ cf.~(\ref{lambda_v}) for the exact definition.

{\bf About the proofs of the Shannon-McMillan theorems.}
We will proceed in several steps, first for lines, for polygons and finally for curves.

{\it (i) Lines:}
The lattice approximation of a line with slope $\lambda\in [0,1]$ and 
$y$-intercept $a$ on an interval $[0,n]$ is defined by 
%\begin{equation}\label{LatApprox}
$L_{\lambda,a}(I)=\{(z,[\lambda z+a])\,\vert\, z=1,...,n\}.$
%\end{equation}
We want to prove the convergence of the rescaled information content
$-(n+1)^{-1}\log P\big[\w_{L_{\lambda,a}([0,n])}\big]$
along sucessively larger segments of the line.
To make this problem accessible to ergodic theory, we have
to find a transformation which captures the stair climbing
pattern along the lattice approximation of the line. If the
slope is rational, the steps become periodic, and we proceed
by combining a finite number of different transformations. 
In the case of an irrational slope, no such simplification is possible.
We need to keep track not only of the integer part but 
also of the fractional part $\{ \lambda z+a\}$ in each step. 
We can realize this by introducing the skew product transformation
\begin{eqnarray*}
      S_\lambda:& \T\times\W&\longrightarrow \T\times\W\\
      \label{SkewProduct}
      & (t,\w) &\longmapsto
      \big(\t_{\lambda}(t),\vartheta_{(1,[t+\lambda])}\w\big),
\end{eqnarray*}
where $\T$ is the one-dimensional torus, equipped with the Borel
$\sigma$-algebra and the Haar measure, and $\t_{\lambda}$ is the
translation by $\lambda.$ Using appropriate ergodic theorems on the
product space with the skew product we obtain a Shannon-McMillan 
theorem along the lattice approximation of the line.  In view of the 
extension of this result to polygons, we further need a variation of this result.
Instead of lattice approximations increasing parts of a line, we use
lattice approximations of blow-ups of a line segment.

{\it (ii) Polygons:}
Consider a polygon $\pi,$ parametrized on $[0,T],$ and its
blowups $B_n\pi(t)=n\pi\big(t/n \big),$ parametrized on $[0,nT].$ 
We study the sequence
%\begin{equation}
$-\log P\big[\w_{L_n^{\pi}}\big]$ $(n\in\N)$
%\end{equation}
of rescaled informations of $P$ restricted to its
lattice approximations $L_n^{\pi}.$ 
Conditioning site by site, the problem can be reduced
to the Shannon-McMillan theorem along the edges,
which is essentially covered by the Shannon-McMillan
theorem along lines that we already established.
One difficulty remains, which is getting around the corners.
It can be overcome by the technique which \cite{FoOrt88} 
used in the case of boxes.  We need the additional assumption 
of a strong form of a $0$-$1$ law (Definition \ref{stzeroone}).
Under this condition, the entropy along a polygon is represented
as a mixture of the entropies of its edges (Theorem \ref{ShMcMpolygon}).
%n the last part of section \ref{gs}

{\it (iii) Curves:} 
Our last step is the entropy along a piecewise smooth curve.
After having established the results for lines and polygons,
this part has become easy.
We obtain (\ref{MainResult}) by approximating the curve with polygons. 

{\bf Relative entropy and large deviations.}
Shannon-McMillan theorems for the
{\it specific relative entropy} \ $h(Q, P)$ of two
probability measures $Q$ and $P$ on the sequence
space $\Upsilon^{\{1,2,\dots\}}$
are based on the functions
$-\log (dQ/dP)\,[\,\w_{\{1,\dots,n\}}]\ (n\in\N)$
describing the {\it relative information} gained
from the first $n$ trials of an experiment.
They are a key tool in the
search for estimates in the theory of large deviations.
By a {\it large deviation} we mean a rare event,
or an untypical behavior occuring in a random sequence.
% ohne s
%The idea of a rare event is central for the theory of
%{\it large deviations.}
Consider the {\it empirical distributions}
%\begin{equation*}
  $\mu_n(\w):= n^{-1}\sum_{i=1}^n \delta_{\w_i}$ $(n\in\N)$
%\end{equation*}
of a stationary random sequence $\w_i\ (i\in\N).$
If the measure $P$ is ergodic then
$\mu_n$ converges to the marginal distribution $\mu$ of $P,$
$P$-almost surely and in $\L^1(P).$ Large deviations are
events like $[\mu_n\in A],$ $A$ being a set in the space of
probability measures on $\Upsilon$
whose closure does not contain $\mu.$

The aim of large deviation theory is to find
lower and upper bounds which describe the asymptotic decay of
the probabilities of such large deviations.
In the classical case of a sequence of
independent and identically distributed random variables, the
decay of large deviations of the empirical distribution
is described by Sanov's theorem. Cram{\'e}r's
theorem addresses similar questions for the
empirical averages.
As a third level for investigating large deviations,
Donsker and Varadhan \cite{DoV75} initiated the investigation
of large deviations of empirical processes.

We replace the random sequence by a random
field, and the empirical processes by the {\it empirical fields}
%\begin{equation*}
  $R_n(\w):=\vert V_n\vert^{-1}
  \sum_{i\in V_n}\delta_{\vartheta_i\w}$ $(n\in\N),$
%\end{equation*}
where $\vartheta_i\ (i\in\Z^d)$ denotes the group the shift
transformations.
Comets \cite{Com86}, F{\"o}llmer and Orey \cite{FoO88},
and Olla \cite{Oll88} found the following
large deviation principle for the empirical fields of
a stationary Gibbs measure:
For any open subset $A$ of the space $\M_1(\W)$,
of probability measures on $\W=\Upsilon^{\Z^d},$
\begin{equation}\label{FoOreyUB}
\liminf_{n\to\infty} \frac{1}{\vert V_n \vert}\log P (R_n\in A)
\ge-\inf_{Q\in A\cap M_1(\W)}\ h (Q,P),
\end{equation}
and for any closed set $C\in\M_1(\W),$
\begin{equation}\label{FoOreyLB}   
\limsup_{n\to\infty} \frac{1}{\vert V_n \vert}\log P (R_n\in
C)\le
-\inf_{Q\in C\cap M_1(\W)}\ h (Q,P),
\end{equation}
where the {\it rate function} is based on the specific relative entropy
$
h(Q,P):=\lim_{n\to\infty}
%\frac{1}{
{\mid V_n\mid}^{-1}H_{V_n}(Q,P).
$

{\bf Large deviations in the phase transition regime.} 
In the case of phase transition, there exists more than one Gibbs measure with
respect to the same potential.  Due to the variational principle for Gibbs measures
(cf.~Lanford and Ruelle \cite{LaR69} and, in purely information theoretical terms,
F{\"o}llmer \cite{Foe73}), the specific relative entropy of
$P$ to another stationary Gibbs measure $Q$ with the same
interaction potential vanishes.
Thus, the relative entropy $h(Q,P)$ appearing in
\eqref{FoOreyUB} and \eqref{FoOreyLB} may be zero
even though $Q$ is not contained in the closure of $A.$
This suggests we need a refined analysis of large deviations
in terms of surface-order rather than volume-order entropies.
Assume in fact that the interaction satisfies the
local Markov property. Then
$H_{V}(Q, P)=H_{\partial V}(Q, P)$
for any finite subset $V$ of $\Z^2,$
where $\partial V$ is the boundary of $V,$ i.e., the
set of all sites outside of $V$ which have distance $1$
to $V$ (cf.~the end of Section 2 in \cite{FoOrt88}).
Consequently, this relative entropy is in fact a
surface-order term, and so it should be rescaled
not by the size of the volume $\vert V \vert$ but by the size
of its surface $\vert\partial V\vert.$ This observation was the main
motivation for introducing the concept of surface entropy, and
for proving the corresponding Shannon-McMillan theorem.

In the context of the two-dimensional Ising model, Schonmann \cite{Sch87} 
showed the existence of surface-order upper and lower bounds for the
large deviations of the emperical means.
For attractive models with a totally ordered state space,
F{\"o}llmer and Ort \cite{FoOrt88} found a lower bound
for the large deviations of the empirical field in terms of
the relative surface entropies along boxes
(recalled as Theorem \ref{LowerBdFO} in this paper).
In their detailed analysis of the two-dimensional Ising model,
Dobrushin, Kotecky, and Shlosman \cite{DKS92}
justified on the basis of local interactions,
that the phase-se\-pa\-rat\-ing curve has the form of a
{Wulff shape}. They proved a large deviation principle
with a rate function in terms of the surface tension along
the Wulff shape. Using different methods,
Ioffe (\cite{Iof94} and \cite{Iof95})
was able to extend their result up to the critical temperature.
The appearance of such shapes suggests to revisit the
approach of \cite{FoOrt88} but using the generalized
surface entropies introduced in Section \ref{gs} instead
of the surface entropies based on boxes.

This extension is carried out in the last section of this work, 
for the case of Gibbs measures with attractive interactions
on a two-dimensional lattice.  We study the large deviations
for its empirical measure. The main result of this part of the paper
is  Theorem \ref{LowerBd}, which provides a lower 
bound in terms of the specific relative entropies along curves.
One of the ingredients in the proof is the well known strategy of 
switching to a measure under which the large deviation becomes 
normal behavior, and then applying a Shannon-McMillan theorem.
Making use of the global Markov property,
we pass from densities restricted to the lattice points inside of
a polygon to densities on the lattice approximations of its boundary.
In this context, we prove an appropriate {\it relative} version of the
Shannon-McMillan theorem, in analogy to the results in Section \ref{gs}.
Other major ingredients in the proof are geometrical observations
concerning the interplay of the random field and the lattice approximations 
of curves.  In particular, in Lemma \ref{LovB} we compute the asymptotic ratio 
of the length of a line segment and its lattice approximation.
These quantities merge into a factor in the lower bound in
Theorem \ref{LowerBd} involving the derivative of the curve.
We will further touch on the case when the Markov property is only satisfied 
with respect a boundary that is a {\it contour} in the sense of statistical mechanics.  
%we prove a similar bound (Theorem \ref{LowerBdc}).

%OOOOOOOOOOOOOOOOOOOOOOOOOOOOOOOOOOOOOOO
%\section*{Outline of the paper}
%OOOOOOOOOOOOOOOOOOOOOOOOOOOOOOOOOOOOOO
{\bf Outline of the paper.}
%(make shorter, part of this belongs in the subsections before... \\
The first section reviews some basic notions and properties 
of discrete random fields, information and entropy.
In the second section, we introduce a specific
entropy of a random field along a line.
In Section \ref{gs} we construct the specific entropy
of a random field along a curve proceeding in three steps: line segments, polygons, curves.
Section \ref{gib} recalls some basic notions about Gibbs measures and phase transitions.  
In Section \ref{lob}  we prove the refined large deviations lower bound.
The case of {\it contour} boundaries is briefly discussed in Subsection
\ref{sec:contours} and in Section \ref{lob}.

%UUUUUUUUUUUUUUUUUUUUUUUUUUUUUUUUUUUUUUUUUUUUUUUUUUUUU
\section{Random fields}\label{rf}
%\markboth{Surface entropy} {Random fields}
%UUUUUUUUUUUUUUUUUUUUUUUUUUUUUUUUUUUUUUUUUUUUUUUUUUUUUU

Consider $\Omega:=\Upsilon^{\Z^d},$ where $\Upsilon$ is a finite
set. For any subset $V$ of $\Z^d$ define $\Omega_V:=\Upsilon^{V}.$
Let $\omega_V$ be the projection of
$\omega$ to $V,$ $P_V$ the distribution of $\omega_V$ with
respect to $P,$ and $\F_V:=\sigma(\omega_V)$ the $\sigma$-algebra
generated by this projection. A probability measure $P$ on
$(\W,\F)$ is called a two-dimensional discrete {\it random field.} 
The transformations $(\theta_{v})_{{v\in\Z^{d}}}$ defined
by $\theta_v\omega(u)=\omega(u+v)\ (u\in\Z^d)$ form the group of
transformations on $\W,$ called {\it shift transformations}. We
assume $P$ is {\it stationary,} that is, invariant with respect to
the shift transformations. 
There are different levels of Markov properties for random fields:
when the subset of the lattice which generates the condition
has to be {\it finite,} and when it can be any type of subset
of the lattice. They both involve the {\it boundary}
$\partial V:=\{\,j\in\Z^d\setminus V\,\vert \, \dist_V(j)=1\,\}$ 
of a subset $V$ of the lattice $\Z^d.$

\begin{definition}\label{MP}
A random field $P$ has the {\rm local Markov property} if, for any finite
$V\subset\Z^d$ and for any nonnegative $\F_V$-measurable $\phi,$
$  E[\phi\,\vert\,\F_{\Z^d\setminus V}] = E[\phi\,\vert\,\F_{\partial V}].$
A random field $P$ which fulfills the local Markov property is
called a {\rm Markov field}. 
If the local Markov property holds for all any $V\subset\Z^d,$ 
then $P$ has the {\rm global Markov property.}
\end{definition}

In Section  \ref{gib}, we will introduce the class of Gibbs
measures in terms of interaction potentials. Any Gibbs measure
belonging to a nearest-neighbor potential is a Markov field.
Examples of random fields which have the local Markov
property but not the global Markov property were given by
Weizs{\"a}cker \cite{Wei83} and Israel \cite{Isr86}. 
$P$ is called {\it tail-trivial} if it fulfills a $0$-$1$ law on the {\it tail field}
\begin{equation}\label{Tail2}
 \Tail:=\bigcap_{V\subset\Z^d\ \text{finite}}\hspace{-4mm}\F_{\Z^d\setminus V}
       =\bigcap_{n\in\N}\F_{Z^d\setminus V_n},
\end{equation}
where
%\begin{equation*}
%  \label{V_n}
 $ V_n:=([-n,n]\cap\Z)^d.$
%\end{equation*}
Due to the spatial structure of a random field, tail-triviality is
equivalent (cf.~Proposition 7.9 from \cite{Geo88}) to a
mixing condition called {\it short-range correlations}\,:
\begin{equation*}\label{ShortRangeWh}
  \sup_{A\in\F_{Z^d\setminus V_{n}}}\vert\,
  P(A\cap B)-P(A)P(B)\,\vert
  \convinfty{n}\,0\qquad  \text{for all } B\in\F.
\end{equation*}
The following strong version of a $0$-$1$ law 
was introduced by F{\"o}llmer and Ort \cite{FoOrt88}.
For $J=\emptyset$ it reduces to the classical $0$-$1$ law on $\F.$ 
Remark 3.2.3 from \cite{FoOrt88} shows that it implies the 
global Markov property provided $P$ has the local Markov property.

\begin{definition}\label{stzeroone}
$P$ satisfies the {\rm strong $0$-$1$ law} if for any subset $J$ of
$\Z^d$ the $\sigma$-algebra $\F_J$ coincides modulo $P$ with
the $\sigma$-algebra
\begin{equation*}
  \F_J^*:=\bigcap_{V\subset {\Z^d}\ \text{{\rm finite}}}\hspace{-4mm}
          \F_{J\cup ({\Z^d}\setminus V)}.
\end{equation*}
\end{definition}

Let $V$ and $W$ be subsets of $\Z^d.$
The {\it information} in $\w$ restricted to $V,$
with respect to $P,$ is given by the random variable
$ \I(P_V)(\omega):=-\log P[\omega_V],$
and the {\it information conditioned} on $\F_W$ is defined as
$  \I\big(P_V[\,\cdot\,\vert\F_W](\omega)\big)
  :=-\log P[\omega_V\vert\,\omega_W].$
The {\it entropy} of $P$ restricted to $V$ is defined as 
\begin{eqnarray*}
  H_V(P):= E[\I(P_V)(\omega)]
  =-\sum_{\w\in\Upsilon^{\Z^d}}P[\w_V]\log P[\w_V],
\end{eqnarray*}
and the {\it conditional entropy} of $P$ to $\F_W$ is
$  H_V(P[\,\cdot\,\vert\F_W])
  :=E[-\log P[\omega_V\vert\,\omega_W]]
  =H(P_V[\,\cdot\,\vert\F_W]).$
The {\it specific entropy} $h(P)$ of $P$ is defined as the limit of
$\vert\,V_n\vert^{-1}H(P_{V_n})$ for $n$ to infinity.
Its existence can be proved by a subadditivity property
(cf., for instance, Theorem 15.12 in \cite{Geo88}),
but it also follows from a Shannon-McMillan theorem
by F{\"o}llmer \cite{Foe73} and Thouvenot \cite{Tho72}. 
They showed that the specific entropy for ergodic $P$ is
$E\big[H\big(P_0[\,\cdot\,\vert\P](\omega)\big)\big],$
where $P_0[\,\cdot\,\vert\P](\w)$ is the conditional distribution of the
random field in the origin with respect to the $\sigma$-algebra $\P$ generated by 
all sites which are smaller than the origin with respect to the lexicographical 
ordering on $\Z^d.$ Moreover, they showed that for all stationary $P$
\begin{eqnarray}\label{limI}
  \frac{1}{\vert\,V_n\vert}\,\I(P_{V_n})\convinfty{n}
  E\big[H\big(P_0[\,\cdot\,\vert\P](\omega)\big)
  \big\vert\J\big]\qquad\mbox{in }\L^1(P),
\end{eqnarray}
where $\J$ is the $\sigma$-algebra of all sets 
invariant with respect to the transformations $\theta_v\ (v\in\Z^d).$

A lemma from \cite{FoOrt88} will be used in some of the proofs in this paper.
We recall it here for the reader's convenience.

\begin{lemma}\label{FoOrt}
Consider $\sigma$-algebras $\B_i\subseteq\B_i^*$ $(i\in\N)$
increasing to $\B_\infty,$ respectively decreasing to
$\B_\infty^*,$ and assume that
%\begin{equation}\label{FoOrt1}
  $\B_\infty=\B_\infty^*\text{ mod }P.$
%\end{equation}
Then for any $\phi\in\L^1(P),$
\begin{equation}
  \lim_{i\to\infty}\sup_{\B_i\subseteq\C_i\subseteq\B_i^*}
  \big\|E[\phi\,\vert\,\C_i]-E[\phi\,\vert\,\B_\infty]\big\|_{\L^1(P)}=0.
\end{equation}
\end{lemma}

%UUUUUUUUUUUUUUUUUUUUUUUUUUUUUUUUUUUUUUUUUUUUUUUUUUUUUUUUU
\section{A Shannon-McMillan theorem along lines}\label{sml}
%\markboth{Surface entropy}{A Shannon-McMillan theorem along lines}
%UUUUUUUUUUUUUUUUUUUUUUUUUUUUUUUUUUUUUUUUUUUUUUUUUUUUUUUUU
From now on we will consider the two dimensional case.
The line with slope $\lambda$ and $y$-intercept $a$ is described by the function
\begin{eqnarray}\label{Defl} 
 l_{\lambda,a}(x):=\lambda x +a\qquad (x\in\R).
\end{eqnarray}
Using $[x]$ and $\{x\}$ for the integer and the fractional part
of $x,$ respectively, the two-sided sequences
$\big([l_{\lambda,a}(z)]\big)_{z\in\Z}$ and
$\big(\{l_{\lambda,a}(z)\}\big)_{z\in\Z}$
are the line's integer and fractional parts at the integer
points $z\in\Z.$ In the case when $0\le\lambda\le 1,$ the
{\it lattice approximation} of $l_{\lambda,a}$ on 
$U\subseteq\Z$ and on an interval $I\subseteq\R$ are given by
\begin{equation}\label{DefL1}
  L_{\lambda,a}(U):=
  \big\{(z,[l_{\lambda,a}(z)])\,\vert\,z\in U\big\}
  \qquad\mbox{and}\qquad
  L_{\lambda,a}(I):=  L_{\lambda,a}(I\cap\Z).
\end{equation}
In the case when $-1\le\lambda<0,$ we use the lattice
approximation $L_{\lambda,a}(I):=-L_{-\lambda,a}(I).$ If
$\vert\lambda\vert >1,$ we represent the line as a function of
the $y$-axis with the new slope $1/\lambda$ (or $0$ in
the case of a parallel to the $y$-axis) and proceed as before.

We want to show the $\L^1(P)$-convergence of the sequence of rescaled
information along successively increased parts of the lattice approximation 
to the line
$$
{\vert L_{\lambda,a}([0,n])\vert}^{-1}\,
\I \big(P_{L_{\lambda,a}([0,n])}\big)\qquad (n\in \N).
$$
In order to make our problem accessible to ergodic theory, we need
to create a transformation that follows the stair climbing pattern
along the lattice approximation of the line.  Lemma $\ref{SthetaIt}$
will introduce a transformation that follows the desired path.
To get there we need a mechanism to keep track not only of the
integer part, but also of the fractional part. Let
$\t_{\lambda}(t):=\{t+\lambda\}\ (t\in\T)$
be the translation by $\lambda$ on the torus $\T:=[0,1]$
with its ends identified and equipped with the Borel $\sigma$-algebra
$\B$ and the Lebesgue measure $\mu.$ Consider the product space
$\T\times\Omega,$ equipped with the product $\sigma$-algebra
$\ov \F,$ and the product measure $\ov P=\mu\otimes P.$
We will now show a few technical lemmata that are needed to proove
the main results in this section.

\begin{lemma}\label{laLa}
For $\lambda\in\R, z,\widetilde{z}\in\Z, a\in\T$ and $I\subset\Z$
we have the following:
\vspace{-2mm}
\begin{enumerate}
\item[(i)] $\t_{\lambda}^z(a)=\{l_{\lambda,a}(z)\}.$
\item[(ii)] The function $\t_{\lambda}^{z}$ has a unique zero. 
More explicitly:
            If $z$ and so $\lambda$ are both positive or both negative, the zero is at
            $a=1-\{z\lambda\}.$ If one is negative
            and the other is positive, the zero is at $a=-\{z\lambda\}.$
            If one of them is zero, then the unique zero is at $a=0.$
\item[(iii)] $l_{\lambda,a}(z+\widetilde{z})=l_{\lambda,a}(z)+\lambda\widetilde{z}$
             \quad and \quad
             $l_{\lambda,a+\widetilde{z}}(z)=l_{\lambda,a}(z)+\widetilde{z}.$
\item[(iv)] $[l_{\lambda,a}(z+\widetilde{z})]=[l_{\lambda,a}(z)]+
             [l_{\lambda,\t_{\lambda}^z(a)}(\widetilde{z})]$
             \quad and \quad
            $[l_{\lambda,a+\widetilde{z}}(z)]=[l_{\lambda,a}(z)]+\widetilde{z}.$
\item[(v)] $L_{\lambda,a}(I+z)=L_{\lambda,\t_{\lambda}^z(a)}(I)+L_{\lambda,a}(z).$
\item[(vi)] $L_{\lambda,a+z}(I)=L_{\lambda,a}(I)+(0,z).$
\end{enumerate}
\end{lemma}
\pf
%\vspace{-3mm}
%\begin{enumerate}
(i) By induction over $z.$
(ii) The case $z=0$ and the case $\lambda=0$ are trivial.  
            Let $z\in\Z\setminus\{0\}.$ By (i), $\t_{\lambda}^z$
            has a zero at $a$ if and only if $\{a+z\lambda\}=0.$
            The latter is equivalent to $a+\{z\lambda\}\in\Z\ (*).$
            If $z$ and $\lambda$ are both positive then
            $0<a+\{z\lambda\}<2.$ 
            Therefore, condition $(*)$ is equivalent to 
            $a=1-\{z\lambda\}.$ If $z$ and $\lambda$ are both
            negative, $\{z\lambda\}$ is again positive and the
            argument works as well. If one is negative and
            the other positive, then
            $-1< a+\{z\lambda\}<1,$ and $(*)$ is
            equivalent to $a+\{z\lambda\}=0.$ 
(iii) $l_{\lambda,a}(z+\widetilde{z})= \lambda(z+\widetilde{z})+a
=l_{\lambda,a}(z)+\lambda\widetilde{z}\quad\mbox{and}\quad
               l_{\lambda,a+\widetilde{z}}(z)=
               \lambda z +a +\widetilde{z}=l_{\lambda,a}(z)+\widetilde{z}.$
(iv) $[l_{\lambda,a}(z+\widetilde{z})]
   =\big[[l_{\lambda,a}(z)]+\{l_{\lambda,a}(z)\}+\lambda\widetilde{z}\big]
   =[l_{\lambda,a}(z)]+[\t_{\lambda}^z(a)+\lambda\widetilde{z}]
   =[l_{\lambda,a}(z)]+[l_{\lambda,\t_{\lambda}^z(a)}(\widetilde{z})],$
   using (iii), (i) and (\ref{Defl}). The second statement is an immediate consequence of (iii).
(v) Using (\ref{DefL1}) and (iv), $L_{\lambda,a}(I+z)=
    \big\{(\widetilde{z},[l_{\lambda,a}(\widetilde{z})])
    \big\vert\,\widetilde{z}\in I+z\big\}
    =\big\{(\widetilde{z}+z,[l_{\lambda,a}(\widetilde{z}+z)])\big\vert
    \,\widetilde{z}\in I\big\}
    =\big\{(\widetilde{z},[l_{\lambda,\t_{\lambda}^z(a)}(\widetilde{z})])+
    (z,[l_{\lambda,a}(z)])\big\vert\,\widetilde{z}\in I\big\}
    = L_{\lambda,\t_{\lambda}^z(a)}(I)+L_{\lambda,a}(z).$
(vi) $L_{\lambda,a+z}(I) =\big\{(\widetilde{z},
   [l_{\lambda,a+z}(\widetilde{z}))\,\vert\,
    \widetilde{z}\in I\big\}
    =\big\{(\widetilde{z},
   [l_{\lambda,a+z}(\widetilde{z}))+(0,z)\,\vert\,
    \widetilde{z}\in I\big\}
    = L_{\lambda,a}(I)+(0,z),$
   using (\ref{DefL1}) and (iv).
%\end{enumerate}
\qed

\begin{lemma}\label{SthetaIt}
The iterates of the transformation
$S_{\lambda}(a,\omega):=(\t_{\lambda}(a),\theta_{(1,[a+\lambda])}\omega)
\ (a\in \T,\w\in\W)$
are given by
$S_{\lambda}^n(a,\omega)=\big(\t_{\lambda}^n(a),\theta_{L_{\lambda,a}(n)}\omega\big)$
for all $n\in\N_0.$
\end{lemma}
\pf With
$\kappa(a):=(1,[a+\lambda]),$ we get
$S_{\lambda}(a,\omega)=(\t_{\lambda}(a),\theta_{\kappa(a)}\omega)$ and
$S_{\lambda}^n(a,\omega)=(\t_{\lambda}^n(a),\theta_{\kappa_n(a)}\omega)$
where $\kappa_n;=\ErgSum{\kappa}{\t_{\lambda}}$ $(n\in\N).$  
It remains to show that
$\kappa_n(a)=(n,[l_{\lambda,a}(n)])$ for all $a\in\T.$
For the first component this is obvious. For the second component it follows by
induction: Trivial for $n=0.$ Then, because of Lemma \ref{laLa}(iv),
$  \kappa_{n+1}^{(2)}(a)
  =\kappa_{n}^{(2)}(a)+\kappa^{(2)}(\t_{\lambda}^n(a))
  =[l_{\lambda,a}(n)]+[\t_{\lambda}^n(a)+\lambda]
  =[l_{\lambda,a}(n)]+[l_{\lambda,\t_{\lambda}^n(a)}]=[l_{\lambda,a}(n+1)].$  
\qed

We use the short forms
%\begin{equation}
    $\P_{\lambda,a}:=\F_{L_{\lambda,a}((-\infty,-1])}$ and $\P_{\lambda,a}^{(i)}:=\F_{L_{\lambda,a}((-i,-1])}$ $(i\in\N)$
%\end{equation}
for the various {\it past} $\sigma$-algebren indexed by the lattice
approximation along of line with slope $\lambda$ and intercept $a.$
They play a central role in the representation
of the limits in the following two Shannon-McMillan type theorems.
%To avoid triple indices we will sometimes skip $\lambda$ in the index of $\t.$
%and write $L(t,I)$ in place of $L_{\lambda,t}(I).$
We further define the functions
\begin{equation}\label{F_i1}
F(a,\w):=
\I\big(P_{0}[\,\cdot\,\vert\P_{\lambda,a}]\big)(\omega), \quad
F_i(a,\w):=
\I\big(P_{0}[\,\cdot\,\vert\P_{\lambda,a}^{(i)}]\big)(\omega) \,\,\,(i\in\N)
\quad\mbox{with } a\in\T, \w\in\W.
\end{equation}

\begin{lemma}\label{FiSi}
For all $a\in\T, \omega\in\Omega,$ and $n\in\N,$\ \
$
 \I(P_{L_{\lambda,a}([0,n])})(\omega)=
%\sum_{i=1}^n F_i\big(\t^i(a),\theta_{(i,[l_{\lambda,a}(i)])}\omega\big)=
  \sum_{i=0}^n F_i\circ S_\lambda^i(a,\omega).
$
\end{lemma}

\pf By conditioning and shifting 
$$
P[\omega_{L_{\lambda,a}([0,n])}]=\prod_{i=0}^n
P[\omega_{L_{\lambda,a}(i)}\vert\omega_{L_{\lambda,a}([0,i-1])}]
= \prod_{i=0}^n P[\omega_{(0,0)}
 \vert\omega_{L_{\lambda,a}([0,i-1])-L_{\lambda,a}(i)}]
 \circ\theta_{L_{\lambda,a}(i)}.
$$
By Lemma \ref{laLa}(v), 
$L_{\lambda,a}([0,i-1])-L_{\lambda,a}(i)=L_{\lambda,\t_\lambda^i(a)}([-i,-1]),$
so 
\begin{eqnarray*}%\label{IasAv}
  \I(P_{L_\lambda,a}([0,...,n])(\omega)=\sum_{i=0}^n
  \I\big(P_{0}[\,\cdot\,\vert\P_{\lambda,\t_\lambda^i(a)}^{(i)}]\big)
  (\theta_{L_{\lambda,a}(i)}\omega),
\end{eqnarray*}
and Lemma \ref{SthetaIt} concludes the proof.
\qed

%To apply ergodic theorems to the right-hand side in Lemma \ref{FiSi},
We will also need the following result about the asymptotic behavior of the functions $F_i$
%$\I\big(P_0[\,\cdot\,\vert\F_{L(\t^i(a),[-1,-i])}]\big)$ 
as $i$ goes to infinity and about their properties as functions of the first parameter. 
Recall that $P_0$ is the marginal distribution of $P$ restricted to the origin.

\begin{lemma}\label{F_icontin}
Assume that
\begin{equation}\label{F_icontinCond}
\forall\, A\in\F_{\Z^2\setminus \{(0,0)\}}:\quad
P_0[\,\cdot\,\vert A]>0.
%\quad\mbox{for}\ P\mbox{-almost all}\ \w\in\W.
\end{equation}
Then, $(F_i)_{i\in\N}$ converge to $F,$ $P$-almost surely, in $\L^1(P)$ in $\w\in\W.$
For any fixed $\w\in\W,$ the functions
$F_i(\cdot,\w)$ $(i\in\N)$ are piecewise constant on $\T.$ 
If $\lambda$ is rational then $F(\cdot,\w)$ is piecewise constant as well,
and the convergence is uniform in $t.$
If $\lambda$ is irrational, if $P$ fulfills the strong $0$-$1$ law and if
\begin{equation}\label{condPbdd}
\exists\, c>0\ \ \forall\, A\in\F_{\Z^2\setminus \{(0,0)\}}:\quad
P_0[\,\cdot\,\vert A]>c,
%\quad\mbox{for}\ P\mbox{-almost all}\ \w\in\W,
\end{equation}
then $F(\cdot\,,\w)$ is Riemann-integrable on $\T,$ and the convergence is
uniform in $t.$

More precisely, the set of points where the function 
$F_i(\cdot,\w)$ may be discontinuous is given by
$\big\{t_{\lambda,\nu}\big\vert\,\nu=-1,...,-i\big\},$
where $t_{\lambda,\nu}$ is the unique zero
of  $\t_{\lambda,\nu}$ on $\T$ (cf.~Lemma \ref{laLa}(ii)).
For rational $\lambda,$ the set of potential points 
of discontinuities of $F_i(\cdot,\w)$ and $F(\cdot,\w)$ is
$\big\{ t_{\lambda,\nu}\big\vert\,\nu=-1,...,-(q\wedge i)\big\},$
where $p/q$ is the unique representation of $\lambda$
with integers $p\in\Z$ and $q\in\N$ having no common divisor.
\end{lemma}

\noindent {\bf Proof of the Lemma.}
Fix any $t\in\T$. Since the $\sigma$-algebras
$\big(\F_{L_{\lambda,t}([-1,-i])}\big)_{i\in\N}$
form an increasing family,
$ \Big(P \big[\w_{(0,0)}\,\big\vert\,
  \w_{L_{\lambda,t}([-1,-i])}\big]\Big)_{i\in\N}$
is a martingale, so that we obtain by the convergence theorem
for martingales that
\begin{equation*}
  P \big[\w_{(0,0)}\,\big\vert\,
  \w_{L_{\lambda,t}([-1,-i])}\big]\convinfty{i}
  P \big[\w_{(0,0)}\,\big\vert\,
  \w_{L_{\lambda,t}([-1,-\infty))}\big]
\end{equation*}
$P$-almost surely and in $\L^1(P).$
By (\ref{F_icontinCond}) this remains true when we take logarithms
on both sides, which proves the first assertion of the lemma.

Now fix $\w\in\W$ and $t,\widetilde t\in\T,$ and find a (sufficient) condition 
under which $F_i(t,\w)=F_i(\widetilde t,\w).$
The only influence that the variable $t$ actually has on $F_i$ 
is its effect on the set $L_{\lambda,t}([-1,-i])$ of sites
we condition on.  By (\ref{DefL1}),
$L_{\lambda,t}\,([-1,-i])= L_{\lambda,\widetilde t}\,([-1,-i])$ if and only if 
\begin{equation}\label{F_icontin_cond}
[l_{\lambda,t}(\nu)]=[l_{\lambda,\widetilde t}(\nu)]\qquad\mbox{for all }
\nu\in\{-1,...,-i\}.
\end{equation}
By Lemma \ref{laLa}(i), $[l_{\lambda,t}(\nu)] =-\lambda \nu-\t_{\lambda}^{\nu}(t),$
so the equality in (\ref{F_icontin_cond}) is equivalent to
$  t-\widetilde t = \t_{\lambda}^{\nu}(t)-{\t}_{\lambda}^{\nu}(\widetilde t).$
This is fulfilled if and only if $t$ and $\widetilde t$
are both either smaller than $t_{\lambda,\nu}$ or both larger than $t_{\lambda,\nu}.$ 
Applying this argument to all $\nu\in\{-1,...,-i\},$ we see that the function
$F_i(\,\cdot\,,\w)$ is piecewise constant, and the set of possible jumps 
is given by $D_i=\big\{t_{\lambda,\nu}\big\vert\,\nu=-1,...,-i\big\}.$

If $\lambda$ is rational, these sets actually become independent of $i,$
for $i$ large enough. With the unique representation $\lambda=p/q$ used above,
and the periodicity of the sequence $(t_{\lambda,\nu})_{\nu\in\N},$ we obtain
$D_i=\big\{t_{\lambda,\nu}\big\vert\,\nu=-1,...,-(q\wedge i)\big\}.$
In particular, the convergence is uniform and the limit $F$ is piecewise constant
in the first variable.

Now assume that $\lambda$ is irrational.
To show that $F$ is Riemann-integrable it suffices to show that the set of 
continuity points has full measure.
We will prove that $F(\cdot,\w)$ is continuous on
$\T\setminus D_\infty,$ where
$D_\infty:=\big\{t_{\lambda,\nu}\vert\,\nu=-1,-2,...\big\}.$
Fix $t_0\in\T\setminus D_\infty$ and let be $\varepsilon>0.$
We apply Lemma \ref{FoOrt} 
with $\B_k=L_{\lambda,t}([-k,-1]),$
and $\B_k^*=L_{\lambda,t}([-k,-1])\cup(\Z^2\setminus V_k),$
where $V_k=[-k,k]^2.$ This gives us a $k_0\in\N$ such that for all
$k\ge k_0$ and $t\in\T$ with
$L_{\lambda,t_0}([-k,-1])=L_{\lambda,t}([-k,-1])\ (**)$ we obtain
\begin{equation*}
  \Big\vert\,
    P\big[\w_{(0,0)}\big\vert\,\w_{L_{\lambda,t_0}([-k,-1])}\big]
  - P\big[\w_{(0,0)}\big\vert\,\w_{L_{\lambda,t}([-k,-1])}\big]
  \,\Big\vert<\varepsilon.
\end{equation*}
By definition of $t_0,$
$\delta:=\min\big\{\vert\, t_0-t_{\lambda,\nu}\vert
\big\vert\,\nu=-1,...,-k\big\}$ is larger than $0,$
and by Lemma \ref{laLa}(i), $(**)$ is true
for all $t\in\T$ for which $\vert t-t_0\vert <\delta.$
So, $P\big[\w_{(0,0)}\big\vert\,\w_{L_{\lambda,t}([-k,-1])}\big]$ 
is continuous in $t=t_0.$
Note that, for all $i\in\N,$ $D_i\subset D_{i+1},$ and
that the maximal height of new jumps added in step $i$,
that is, the steps in $D_{i+1}\setminus D_i,$ decreases with $i.$
%$\mbox{max}\,\Psi(F_{i+1})(D_{i+1}\setminus D_i)\le \mbox{max}\,\Psi(F_i)(D_i).$
This implies the uniformity of the convergence.
Finally, by taking logarithms and using (\ref{condPbdd}),
the above statements are also true for the sequence $(F_i(\cdot\,,\w))_{i\in\N}.$
\qed

%..................\t periodisch..............................
\begin{thm}\label{SMQ}
Let $\lambda$ be rational and $p/q$ its unique representation
with $p\in\Z$ and $q\in\N$ having no common divisor.
Assume that $P$ fulfills condition \eqref{F_icontinCond}.
Then for all $a\in\R,$
\begin{eqnarray}\label{SMQ1}
\frac{1}{n+1}\,\I(P_{L_{\lambda,a}([0,n])})\convinfty{n}
\frac{1}{q}\sum\limits_{\nu=0}^{q-1}E\Big[H\big(P_0\big[
\,\cdot\,\big\vert
\P_{p/q,\,\nu p/q}\big](\omega)\big)\Big\vert\J_{q,p}\Big]
\end{eqnarray}
$P$-almost surely and in $\L^1(P),$ where $\J_{q,p}$ is the
$\sigma$-algebra of all sets which are invariant with respect to $\theta_{(q,p)}.$
In particular, if $P$ is ergodic with respect to $\theta_{(q, p)}$ then the limit equals
$$
\frac{1}{q}\sum_{\nu=0}^{q-1}
E\Big[H\big(P_0\big[\,\cdot\,\vert\P_{p/q,\,\nu p/q}\big](\omega)\big)\Big].
$$
\end{thm}

\pf
% reduce to F via Maker
By Lemma \ref{FiSi}, the left-hand side of \eqref{SMQ1} equals
$  1/(n+1)\sum_{i=0}^n F_i\circ S_\lambda^i(a,\omega).$
Because of  Lemma \ref{F_icontin} and Maker's version of the ergodic theorem
(cf.~Theorem 7.4 in Chapter 1 of \cite{Kre85}), it suffices to show that
 \begin{equation}\label{SMQ2} 
   \frac{1}{n+1}\sum_{i=0}^n F\circ S_\lambda^i(a,\omega) 
\end{equation}
converges to the limit in \eqref{SMQ1}.
% proof for F
For each $n\in\N$ there are unique $m\in\N$ and $\eta\in\{0,1,...,q-1\}$ such that $n=mq+\eta.$
The last term can be rewritten as %\eqref{SMQ2}
\begin{displaymath}
  %\ErgAv{F}{S}=
  \frac{m}{n+1}\sum_{\nu=0}^{q-1}
  \frac{1}{m}\sum_{j=0}^{m-1}F\circ S_\lambda^{jq+\nu}
  +\frac{1}{n+1}\sum_{\nu=0}^{\eta}F\circ S_\lambda^{mq+\nu}.
\end{displaymath}
The second addend converges to $0$ as $n$ (and therefore also $m$) goes to infinity.
The first factor converges to $1/q,$ so it remains to study the asymptotic behavior of 
\begin{eqnarray*}
   %\label{Reterg2}
   A_m^{(\nu)}F(a,\w)
    := \frac{1}{m}\sum_{j=0}^{m-1}F\circ S_\lambda^{jq+\eta}(a,\w).
\end{eqnarray*}
Use $\kappa_n$ $(n\in\N)$ defined as in the proof of Lemma \ref{SthetaIt}.
Since $\t_\lambda^q=\mbox{Id}$ we obtain $\kappa_{jq+\nu}=j\kappa_q+\kappa_\nu,$ and
$S_\lambda^{jq+\nu}=\big(\t_\lambda^{\nu},\theta_{\kappa_\nu}\circ(\theta_{\kappa_q})^j\big).$
This yields
\begin{eqnarray*}
     A_m^{(\nu)}F(a,\w)
     =\frac{1}{m}\sum_{j=0}^{m-1}
     F\big(\t^{\nu}(a), \theta^{\kappa_{\nu}(a)}\circ(\theta^{\kappa_q(a)})^{\,j}\w\big).
\end{eqnarray*} 
For $a$ fixed, applying Birkhoff's ergodic theorem to the function 
$f^{(\nu)}_a(\w):= F\big(\t^\nu(a),\theta_{\kappa_\nu(a)}\w\big)$  yields
\begin{displaymath}
  \lim_{m\to\infty} A_m^{(\nu)}F(a,\cdot)
  =E[f^{(\nu)}_t\vert\J_{(q,p)}]
  =E\Big[H\big(P_0\big[\,\cdot\,\big\vert\P_{p/q,\,\nu p/q}\big](\omega)\big)\Big\vert\J_{(q,p)}\Big]
%  =E\big[F\big(\t^\nu(a),\theta_{\kappa_{\nu}(a)}(\cdot)\big)\vert\J_{(q,p)}]
  \quad\mbox{$P$-a.s.~and in $\L^1(P).$}
\end{displaymath}
In the ergodic case $\J_{(q,p)}$ is trivial.
Using the invariance of $P$ under $\theta,$ the last expression reduces to 
$E\big[F\big(\t^\nu(a),\,\cdot\,\big)].$ 
\qed

For the case of an irrational slope $\lambda$ we can show the corresponding result
provided $P$ fulfills some additional conditions, in particular 
the strong $0$-$1$ law (cf.~Definition \ref{stzeroone}).

%....................\t ergodisch..............................
\begin{thm}\label{SMnotQ}
Let $\lambda$ be irrational. Assume that
$P$ fulfills condition 
\eqref{condPbdd}
and the strong $0$-$1$ law.
Then for all $a\in\R,$
\begin{eqnarray}\label{SMnotQLim}
\frac{1}{n+1}\,\I(P_{L_{\lambda,a}([0,n])})\convinfty{n}
\int_0^1 E\big[H\big(P_0[\,\cdot\,\vert
\P_{\lambda,t}](\w)\big)\big]\,dt\qquad\mbox{in}\ \L^1(P).
\end{eqnarray}
\end{thm}

\pf
As in the proof of Theorem \ref{SMQ} the left-hand side of \eqref{SMQ1} can be
rewritten as an ergodic average involving the skew product $S_\lambda.$
Using Lemma \ref{F_icontin} and Maker's version of the ergodic theorem
(cf.~Theorem 7.4 in Chapter 1 of \cite{Kre85}), it suffices to show that for all $i\in\N,$
\begin{eqnarray}\label{SMnotQLim_k}
1/(n+1)\sum_{i=0}^n F_k\circ S_\lambda^i(a,\omega)
\convinfty{n}
\int_0^1 E\big[H\big(P_0[\,\cdot\,\vert\P_{\lambda,t}^{(i)}](\w)\big)\big]\,dt
\qquad\mbox{in}\ \L^1(P)
\end{eqnarray}
For $i\in\N$ fixed, $F_i$ is piecewise constant in $t,$ so there is an $R\in\N,$ 
intervals $U_1,...,U_R$ and measurable functions $f_1,...,f_R$ on $(\W,\F)$ such that
%\begin{equation}\label{F_factors}
$F_i(t,\w)=\sum_{r=1}^R 1_{U_r}\,f_r(\w).$
%\end{equation}
By approximation, this can further be reduced to the case of indicator functions
at the place of the $f_r$'s. It remains to show the convergence for functions of the form
$1_U(t)\,1_A(\w)$ for an interval $U\subset\T$ and a set $A\in\F.$ 

Let $F_i(t,\w)=1_U(t)\,1_A(\w).$ We will show the convergence in $\L^2(P).$ Using 
$$ 
\int_\W 1_A( \theta_{\kappa_i(t)}(\w) 1_A(\theta_{\kappa_j(t)}\w)\,P(d\w) =
%P\big(\theta_{\kappa_i(t)}^{-1}A\cap\theta_{\kappa_j(t)}^{-1}A\big) =
P\big(\theta_{\kappa_i(t)-\kappa_j(t)}^{-1}A\cap A\big),
$$
we see that
\begin{align}
%\begin{split}
 \bigg\| & \frac{1}{n+1} \sum_{i=0}^{n} 
    F\big(\t^i(t), \theta_{\kappa_i(t)}\,\cdot\big) 
   -  \int_0^1 E[F]\,dt \bigg\|_{\L^2(P)}^2 \nonumber\\
  & =  \frac{1}{(n+1)^2}\sum_{i=0}^{n}1_U(\t^i(t))1_U(\t^j(t))
    P\big(\theta_{\kappa_i(t)-\kappa_j(t)}^{-1}A\cap A\big)
   + \  \mu(U)^2 P(A)^2 \label{BHprdlem1} \\
    &\quad -  \frac{2}{(n+1)^2}\sum_{i=0}^{n}\bigg[  
        1_U(\t^i(t))
        \int_{\W} 1_A(\theta_{\kappa_i(t)}\w)\,P(d\w)
         \bigg]
       \mu(U)P(A) \label{BHprdlem1a} 
\end{align}

The next step is to show that the first addend 
in line (\ref{BHprdlem1}) may be replaced by
\begin{equation}\label{BHprdlem2}
 \bigg( \frac{1}{(n+1)}\sum_{i=0}^{n}1_U(\t^i(t))\bigg)^2P(A)^2
\end{equation}
without affecting the asymptotic behavior of (\ref{BHprdlem1}).
% for $n$ going to infinity: 
We may bound
\begin{align*}
  \bigg\vert&\frac{1}{(n+1)^2}\sum_{i,j=0}^{n}1_U(\t^i(t))1_U(\t^j(t))
   \Big(P\big(\theta_{\kappa_i(t)-\kappa_j(t)}^{-1}A\cap A\big)
    -P(A)^2\Big)\bigg\vert\nn\\
  &\le  \frac{1}{(n+1)^2}\sum_{i,j=0}^{n}
    \Big\vert P\big(\theta_{\kappa_i(t)-\kappa_j(t)}^{-1}A\cap A\big)
    -P(A)^2\Big\vert.
\end{align*}
Fix $\varepsilon>0.$ 
Let $\|\,\cdot\,\|$ denote the maximum norm on $\Z^2.$
As a special case of the strong $0$-$1$ law, $P$ fulfills a $0$-$1$ law on the tail field.
By (\ref{ShortRangeWh}) there is an $m\in\N$ such that
\begin{equation*}
  \Big| P\big(\theta_k^{-1}A\cap A\big) - P(A)^2\Big|
  <\frac{\varepsilon}{2}\quad\quad\text{for all }k\in\Z^2
  \text{ with }\|k\|> m.
\end{equation*}
Note that $\| \kappa_i(t)-\kappa_j(t)\|\ge | i-j |$ for all $t\in\T.$ Since
$\lim_{n\to\infty}n^{-2}\big\vert\big\{1\le i,j\le n\,\big\vert\,| i-j |\le m\big\}\big\vert=0$
for all $m\in\N,$ 
we can find an $n_0\in\N$ such that
\begin{equation*}
  \frac{1}{(n+1)^2}\,\sup_{t\in\T}\,\big|\big\{1\le i,j\le
  n\,\big\vert\,\|\kappa_i(t)-\kappa_j(t)\|\le m\big\}\big|
  <\frac{\varepsilon}{2}
  \quad\text{ for all }n\ge n_0.
\end{equation*}
This demonstrates that hat the difference created by the by replacing the 
first addend in \eqref{BHprdlem1} by \eqref{BHprdlem2}
converges to $0$ uniformly with respect to $t.$ 

Since $U$ is an interval, the function $1_U$ is integrable in the sense of Riemann.
Applying a variant of Weyl's theorem for Riemann integrable functions (cf., for instance, 
Theorem 2.6  in Chapter 1 of \cite{Kre85}) implies that, for $n$ to infinity, 
$1/(n+1) \sum_{i=0}^{n} 1_U(\t^i(t))$ converges to $\mu(U)$ uniformly in $t.$ 
So, asymptotically, the sum of the two addends in line \eqref{BHprdlem1}
equals $2\,\mu(U)^2 P(A)^2.$ 

The expression in \eqref{BHprdlem1a} can be simplified to
   $-  2/(n+1)\sum_{i=0}^{n} 1_U(\t^i(t))\mu(U)P(A)^2.$
Applying a variant of Weyl's theorem for Riemann integrable functions again,
this expression converges to $-2\,\mu(U)^2 P(A)^2,$ which concludes the proof. 
\qed

%UUUUUUUUUUUUUUUUUUUUUUUUUUUUUUUUUUUUUUUUUUUUUUUU
\section{A Shannon-McMillan theorem along general shapes}\label{gs}
%\markboth{Surface entropy}{General shapes}
%UUUUUUUUUUUUUUUUUUUUUUUUUUUUUUUUUUUUUUUUUUUUUUUU

Let $P$ be a stationary random field that satisfies the
strong $0$-$1$ law and the condition \eqref{condPbdd}.
The goal of this section is a Shannon-McMillan theorem for a stochastic field 
along the lattice approximations of the blowups of a curve $c$, and an explicit
formula for the limit $h_c(P)$, the {\it specific entropy} of $P$ {\it along} $c.$
This will be done in three steps: for linear segments, for polygons and for curves.
In Subsection \ref{sec:contours}, 
we will further introduce lattice approximations that are {\it contours} in the sense 
of statistical mechanics and sketch corresponding results for them.

Let $c=\big(c^{(1)},c^{(2)}\big)$ be a piecewise differentiable planar curve 
parametrized by $t\in [0,T].$
Assume that the trace of $c$ does not contain the origin and that it hits the $y$-axis in $t=0.$
Let $c^{\prime}$ denote the right derivative of $c.$  The blowups of the curve $c$ are given by
\begin{equation}\label{BlowingupDef}
  B_{\eta} c:[0,\eta T)\longrightarrow \R^2, \quad\quad
  B_{\eta} c(t)=\eta\,c(t/\eta)\quad\quad(\eta>0).
\end{equation}
It will follow from the construction in Subsection \ref{sec:poly} that is enough to consider 
curves that are described by the graph of a function $\phi$ on a segment of one
of the axes. Suppose that $\phi$ is a function on the interval
$[x,\widetilde{x}]$ of the $x$-axis. (The case of the $y$-axis can be
treated analogously.) More precisely, $x=c^{(1)}(0)$ and
$\widetilde{x}=c^{(1)}(T).$ The interval $[x,\widetilde x]$ contains a
finite number $u$ of integers $z,...,z+u,$ where
$u=[x-\widetilde x]$ or $u=[x-\widetilde x]-1.$
In the same way, the blowups $B_{n}c$ can be
represented as graphs of functions $\phi_n$ of
intervals $[x_n,\widetilde x_n].$ 
We obtain by (\ref{BlowingupDef}),
       $ x_n = B_nc^{(1)}(0) =nx=nc^{(1)}(0)$ and
$ \widetilde x_n = B_nc^{(1)}(nT)=n\widetilde x=nc^{(1)}(T).$
Again, the interval $[x_n,\widetilde x_n]$ contains a finite
number $u_n$ of integers $z_n,z_n+1,...,z_n+u_n,$ where
\begin{equation}\label{u_n}
  u_n=[n(\widetilde x-x)]\quad\quad\text{or}\quad\quad
  u_n=[n(\widetilde x-x)]-1.
\end{equation}
In particular, the sequence $(u_n)_{n\in\N}$ goes to infinity.

%UUUUUUUUUUUUUUUUUUUUUUUUUUUUUUUUUUUUUUUUUUUUUUUUUUUUUUUUUUUUUUU
\subsection{Line segments}\label{sec:lines}
%\markboth{Surface entropy}{General shapes}
%UUUUUUUUUUUUUUUUUUUUUUUUUUUUUUUUUUUUUUUUUUUUUUUUUUUUUUUUUUUUUUU

The main result of this section %(cf.~\ref{ShMcMStretchL})
is the convergence of the sequence of renormalized information functions
\begin{equation}\label{InfStretch}
       \frac{1}{\vert L_n(a)\vert}\,\I\big(P_{L_n(a)}\big)(\w)
\end{equation}
along the lattice approximation converges to the entropy $h_{\lambda}(P)$ of $P$ along a line with 
slope $\lambda.$  If $c$ is a line segment, the functions $\phi_n$ are of the form
%\begin{align*}
  $\phi_n(x)= \lambda x+a_n \
  %& \quad\text{for}\quad 
  (x\in[nc^{(1)}(0),nc^{(1)}(T)])$
%  \text{where}& \quad
 where
 $  \lambda=({c^{(2)}(T)-c^{(2)}(0)})/({c^{(1)}(T)-c^{(1)}(0)})$
  %\quad\text{and}\quad 
  and $a_n=n(c^{(2)}(0)-\lambda c^{(1)}(0)).$
%\end{align*}
% \\ @MMM  definition von
Assume $0\le\lambda\le 1.$ As explained at the beginning of Section \ref{sml},
the other cases can be reduced to this case.  At first sight it seems we could 
just apply the results for the specific entropy along a line from the last section. 
But the blowups of the line segments move in space, which has the following consequences:
\vspace{-2mm}
\begin{itemize}
\item[(i)] There is a sequence
      $(a_n)_{n\in\N}$ instead of a constant $a.$
\item[(ii)] The sequence is real-valued, as opposed to the constant
            having values in $\T.$
\item[(iii)] The positions of the lattice points in each step are more difficult to describe. 
      In the case of the line we simply looked at approximating points with
      $x$-values between $0$ and $n.$ Now, $x$-values of the lattice points
      lie in the interval between $z_n$ and $z_n+u_n.$
\end{itemize}
\vspace{-2mm}
The last problem forces us to apply, at each step $n,$ an
additional shift to $\w$ which brings the line segment close to
the origin. These shifts do not affect the $\L^1(P)$-convergence
since the limit is shift invariant. The number of points in
the $n$th step is given by $(u_n)_{n\in\N},$ instead of simply
$n+1$ as in the last section, but this is irrelevant as long as
the sequence goes to infinity. The second problem requires another
shift in each step $n.$ The first point is the most delicate. It
is here that we need the convergence of the ergodic
averages in {\it all} $t,$ rather than just {\it almost} all $t.$ 
%which was discussed in Section 4 of \cite{Bre01}.

\begin{thm}\label{ShMcMStretchL}
In $\L^1(P)$ and uniformly in $a\in\R,$
\begin{equation*}
   %\lim_{n\to\infty}
   \frac{1}{\vert L_n(a)\vert}\,\I(P_{L_n(a)})\convinfty{n} h_{\lambda}(P).
\end{equation*}
\end{thm}

\pf
We start with some technical preparations.
Set $a:=a_1.$ Using the notation for $[x,\widetilde x]$
described above \eqref{u_n}, define for $n\in\N,$
%\begin{equation}\label{DefL}  
    $L_n(a):= L_{\lambda,a_n}(z_n,\dots,z_n+u_n).$
%\end{equation}
The total number of sites in $L_n(t)$ is $u_n+1.$ 
To transform \eqref{InfStretch} into some sort of ergodic average
we first condition on successively smaller parts of $L_n(a).$
A new step begins at $i$ if and only if $\t_{\lambda}^{z_n+i-1}(\{a\})\ge 1-\lambda. $
%\begin{lemma}\label{LDiff}
For all $z, i \in\Z,$ and $a\in\R,$
\begin{equation}\label{LDiff}
     L_{\lambda,a}(z+i) - L_{\lambda,a}(z)=L_{\lambda,\t_\lambda^{z}(\{a\})}(i).
\end{equation}
To prove this, first apply the second equation in Lemma \ref{laLa}(v)
with $a=\{a\}$ and $z=[a],$ and then apply the first equation with
$\widetilde{z}=i,$ %and finally the definition (\ref{DefL}) to obtain
\begin{align*}
  L_{\lambda,a}(z+i)
 & = (z+i,[l_{\lambda,\{a\}}(z+i)]) + (0,[a])\\
 &  = (z+i,[l_{\lambda,\{a\}}(z)]
    + [l_{\lambda, \t_\lambda^{z}(\{a\})}(i)])+ (0,[a])\\
 & = \big(z,[l_{\lambda,a}(z)]\big)
    + \big(i,[l_{\lambda,\t_\lambda^{z}(\{a\})}(i)]\big)
    = L_{\lambda,a}(z) + L_{\lambda,\t_\lambda^{z}(\{a\})}(i).
\end{align*}

We calculate the information in (\ref{InfStretch}) by conditioning site by site along 
$L_n(a).$  We use $\w(i)$ instead of $\w_i$ for easier reading.  Shifting $\w$ to the origin, 
applying \eqref{LDiff} and using the functions defined in \eqref{F_i1}  yields
\begin{align}\label{ShMcM1}  
\I( & P_{L_n(a)})(\w) % \nonumber \\ &
=- \sum_{i=0}^{u_n}
      \log P\Big[\w\big(L_{\lambda,a_n}(z_n+i)\big)
     \,\Big\vert\,\w\big(L_{\lambda,a_n}(z_n+i-1,\dots,z_n)\Big]
          \nonumber \\
%\end{align}
%\begin{align}\label{ShMc0a}  %or add1???
&=- \sum_{i=0}^{u_n} 
    \log P\Big[\w(0,0)
    \,\Big\vert\,\w\big(L_{\lambda,a_n}(z_n+i-1,...,z_n)-L_{\lambda,a_n}(z_n+i)\big)\Big]
    \circ \theta_{L_{\lambda,a_n}(z_n+i)} 
    \nonumber \\
&=- \sum_{i=0}^{u_n}   
   \log P\Big[\w(0,0)
   \,\Big\vert\,\w\big(L_{\lambda,\t_\lambda^{z_n+i}(\{a_n\})}(-1,...,-i)\big)\Big]
   \circ\theta_{ L_{\lambda,\t_\lambda^{z_n}(\{a_n\})}(i)}
   \circ\theta_{ L_{\lambda,a_n}(z_n) }\\
&=- \sum_{i=0}^{u_n}   
     F_i\big(\t_\lambda^i(\t_\lambda^{z_n}(\{a_n\})),
      \theta_{L_{\lambda,\t_\lambda^{z_n}(\{a_n\})}(i)}\circ
      \theta_{L_{\lambda,a_n}(z_n)}\,\w\big).
      \nonumber 
\end{align}  
%\end{equation*}
Applying Lemma \ref{SthetaIt},
putting all together back in (\ref{ShMcM1}) and renormalizing yields
\begin{align}\label{ShMcM6}
  \frac{1}{u_n+1}\,\I(P_{L_n(a)})(\w)
 &=  \frac{1}{u_n+1}\sum_{i=0}^{u_n}
 F_i\circ S_\lambda^i\big(\t_\lambda^{z_n}(\{a_n\}),
 \theta_{L_{\lambda,a_n}(z_n)}\,\w\big).
\end{align}
To prove the convergence, we have to distinguish the case when $\lambda$
is rational from the cases when it is irrational, because this
determines whether $\t_{\lambda}$ is periodic or uniquely ergodic. 
We proceed as in the proof of Theorem \ref{SMQ} and Theorem \ref{SMnotQ}, respectively.
\qed

%UUUUUUUUUUUUUUUUUUUUUUUUUUUUUUUUUUUUUUUUUUUUUUUUUUUUUUUUUUUUUUU
\subsection{Polygons}\label{sec:poly}
%\markboth{Surface entropy}{General shapes}
%UUUUUUUUUUUUUUUUUUUUUUUUUUUUUUUUUUUUUUUUUUUUUUUUUUUUUUUUUUUUUUU

The next step is to define the entropy along a polygon, that is a
piecewise linear curve $\pi: [0,T]\rightarrow\R^2.$ 
Assume further that $\pi$ fulfills the other assumptions on $c$
stated at the beginning of this section.  Let $R$ be the number of edges 
of $\pi.$ We can find slopes $\lambda^{(r)}\in (-1,1],$ constants $t^{(r)}\in\R,$ 
and intervals $I^{(r)}$ of the $x$- or the $y$-axis such that
\begin{equation}
  \pi\big([0,T]\big)=\bigcup_{r=1}^R l_{\lambda^{(r)},t^{(r)}}(I^{(r)}),
\end{equation}
with $l_{\lambda,t}$ as defined in (\ref{Defl}) as a function of
the $x$- or of the $y$-axis. Proceeding the same way for the blowups
$B_n\pi\;(n\in\N)$ as defined in (\ref{BlowingupDef}), we choose $t_n^{(r)}\in\R$ and
$I_n^{(r)}\subset\R,$ such that
\begin{equation*}
  B_n\pi\big([0,T]\big)=
  \bigcup_{r=1}^R l_{\lambda^{(r)},t_n^{(r)}}(I_n^{(r)}),
\end{equation*}
The lattice approximations of the edges combine to a lattice approximation of $B_n\pi:$
\begin{equation}\label{LBnpi}
  L_n^\pi:=\bigcup_{r=1}^R L_{\lambda^{(r)},t_n^{(r)}}
\end{equation}

\begin{thm}\label{ShMcMpolygon}
The lattice approximations converge. More precisely, 
\begin{equation}\label{ShMcMpolygon1}
  \frac{1}{\len L_n^\pi}\,\I(P_{L_n^\pi})\convinfty{n}
  \frac{1}{\len\pi}
  \sum_{r=1}^R \len\pi^{(r)}\, h_{\lambda^{(r)}}(P)\qquad\mbox{in }\L^1(P).
\end{equation}
\end{thm}
In some contexts it is more convenient to express the limit as an integral
with respect to $t$ rather than as a sum.  
Let $\pi'(t)$ denote the right derivative of $\pi.$
Then the limit can be written as
\begin{equation}\label{IntReprhpi}
  \frac{1}{\len\pi}\int_0^T h_{\pi'(t)}(P)\,dt.
\end{equation}

\pfofthm
Use $\w(i)$ for $\w_i$ and define the sets
%\begin{equation}\label{Enr}
  $E_n^{(r)}:=L_{\lambda^{(r)},t_n^{(r)}}(I_n^{(r)})$
  $(r\in\{1,\dots,R\}).$
%\end{equation}
By conditioning, 
\begin{equation}\label{pfpoly1}
  \I\big(P_{L_n^\pi}\big)(\w)
  =\sum_{r=1}^R \log P\big[\w(E_n^{(r)})
         \,\Big\vert\,\w(E_n^{(r-1)},\dots,E_n^{(1)})\big].
\end{equation}
Fix $r\in\{1,\dots,R\}.$ Omit the index
$r$ when there is no risk of confusion (for example
$\lambda:=\lambda^{(r)}, t_n:=t_n^{(r)}, E_n:=E_n^{(r)}$) and use the short form
%\begin{equation*}%\label{EnRest}
$\breve{E}_n:=E_n^{(r-1)}\cup\dots\cup E_n^{(1)}$
%\end{equation*}
for the lattice approximations of the edges of the polygon which come prior to $E_n^{(r)}.$ 
We will condition successively on the elements of $E_n^{(r)}.$ Denoting the integers in 
$I_n^{(r)}$ by $z_n, z_n+1,\dots,z_n+u_n$ as in (\ref{u_n}),we obtain for the $r$th addend in (\ref{pfpoly1})
\begin{align*}
        \log P\Big[  \w(E_n) \,\Big\vert\,\w\big(\breve{E}_n\big)\Big]
        =\sum_{i=0}^{u_n}
        \log P\Big[\w\big(L_{\lambda,t_n}(z_n+i)\big)
        \,\Big\vert\,
        \w\big(L_{\lambda,t_n}(z_n+i-1,\dots,z_n),
        \breve{E}_n\big)\Big].
\end{align*}
Shifting by $v_{ni}(t):=L_{\lambda,t_n}(z_n+i)$ yields
\begin{align*}
     \sum_{i=0}^{u_n} &
     \log P\Big[\w(0,0)\Big\vert\w\big(\L_{\lambda,t_n}(z_n+i-1,\dots,z_n)
         -v_{ni}(t),\breve{E}_n- v_{ni}(t)\big)\Big]\circ\theta_{v_{ni}(t)}.
\end{align*}
By \eqref{LDiff}, this equals
%$v_{ni}(t)=L_{\lambda,t_n}(z_n)+L_{\lambda,\t_\lambda^{z_n}(\{t_n\})}(i),$ so we obtain
\begin{align}\label{pfpoly4}  
   \sum_{i=0}^{u_n} &
        \log P\Big[\w(0,0)\,\Big\vert\,
        \w\big(L_{\lambda,\t_{\lambda}^{z_n+i}(\{t_n\}) }
        (-1,\dots,-i),\breve{E}_n-v_{ni}(t)\big)\Big]  
        \nonumber \\
  & \phantom{\log P\Big[\w\big(0,0\big)\,\Big\vert\,\,\Big)
     \w L_{\lambda,\t_{\lambda}^{z_n+i}} (\{t_n\})(-1,\dots,-i) }
         \circ\theta_{L_{\lambda,\t_\lambda^{z_n}(\{t_n\})}(i)}
         \circ\theta_{L_{\lambda,t_n}(z_n)}.
\end{align}
This expression is similar to (\ref{ShMcM1})  %and (\ref{add2})
except for the additional conditionings on the sites $\breve{E}_n-v_{ni}(t).$ 
We will show that these conditions disappear asymptotically. 
The argument will be given in detail for the first summand; 
it is similar for the remaining ones.
Let $\alpha$ be the minimum angle between any neighboring edges of
the polygon $\pi$ and let $d_n$ be the minimum distance between an
edge of the $n$th blowup $B_n\pi$ of the polygon and any of its
nonneighboring edges. Also, let $H_n$ be the hexagon defined as
follows: $H_n$ is symmetric around $E_n^{(r)},$ two sides are
parallel to $E_n^{(r)}$ at a distance $d_n/2.$ The other sides
reach from the endpoint of the first two to the endpoints of
$E_n^{(r)},$ and they intersect at an angle $\alpha$.
Observe that $\breve{E}_n\subset\Z^2\setminus H_n,$ and therefore
$\breve{E}_n-v_{ni}(t)\subset\Z^2\setminus\big(H_n-v_{ni}(t)\big).$
Define the $\sigma$-algebras
$\B_i(t) :=\F\big(L_{\lambda,\t_\lambda^i(\{t\})}(-1,\dots,-i)\big),$ 
$\B_\infty(t) :=\F\big(L_{\lambda,\t_\lambda^i(\{t\})}(-1,-2,\dots)\big)$ and 
$\B_i^*(t) :=\F\big(L_{\lambda,\t_\lambda^i(\{t\})}(-1,\dots,-i)\cup
                   \Z^2\setminus \big(H_n-v_{ni}(t)\big)\big).$
The sequence $\big(\B_i(t)\big)_{i\in\N}$ is increasing to
$\B_\infty(t),$ and the sequence $\big(\B_i^*(t)\big)_{i\in\N}$ is decreasing to 
$\B_\infty^*(t):=\bigcap_{i\in\N}\B_i^*(t). $ By the strong $0$-$1$ law, 
$\B_\infty^*(t)=\B_\infty(t)\ \text{mod}\ P.$ By Lemma \ref{FoOrt}, 
\begin{align*}
\lim_{i\to\infty}
\Big\|
 \log P\Big[\w(0,0)
        \,\Big\vert\, &
        \w\big(L_{\lambda,\t_\lambda^i(t)}(-1,\dots,-i),
        \breve{E}_n-v_{ni}(t)\big)\Big]\\
   & -  \log P\Big[\w(0,0)
        \,\Big\vert\,
        \F\big(L_{\lambda,\t_\lambda^i(\{t\})}(-1,\dots,-i)\big)
        \Big](\w)
 \Big\|_{\L^1(P)}
=0.
\end{align*}
Proceeding with (\ref{pfpoly4}) as  with (\ref{ShMcM1}) and using that, for all $r\in\{1,\dots,R\},$  
$ \len E_n^{(r)}/\len L_n^{\pi}$ asymptotically equals $\len\pi^{(r)}/\len\pi$ concludes the proof. 
\qed

%UUUUUUUUUUUUUUUUUUUUUUUUUUUUUUUUUUUUUUUUUUUUUUUUUUUUUUUUUUUUUUU
\subsection{Curves}\label{sec:curves}
%\markboth{Surface entropy}{General shapes}
%UUUUUUUUUUUUUUUUUUUUUUUUUUUUUUUUUUUUUUUUUUUUUUUUUUUUUUUUUUUUUUU

To use the results from the previous section, we need to relate the derivatives of
the curve with a slopes of lines.    
Let $v\in S^1=\{w\in\R^{2}\,\vert\,\,\vert w\vert=1\},$
and $\alpha$ the angle from the positive $x$-axis to the vector $v.$
If $\vert\alpha\vert\le\pi/4$ or $\vert\alpha\vert\ge 3\pi/4$
then describe the line in the direction of $v$ by a function
of the $x$-axis; otherwise describe it as a function of the $y$-axis.
We assign any $v\in S^1$ a specific entropy
\begin{equation}\label{lambda_v}  
   h_{v}(P):=h_{\lambda(v)}(P),  
   \mbox{where}\   \lambda(v):=\min(\vert tg\,\alpha\vert,\vert ct\,\alpha\vert).
\end{equation}

\begin{thm}\label{ShMcMCurve}
Let  $c:[0,T]\longrightarrow\R^2 $ be a piecewise continuously
differentiable curve.  Assume the trace does not contain the origin.  
Let $\pi_n:[0,nT]\longrightarrow\R^2\ (n\in\N) $ be a sequence of polygons such that
\begin{align}\label{PolyCond}
  \frac{1}{\len\pi_n}\,
  \sup_{t\in [0,nT]}\big\vert(B_n^{-1}\pi_n)'(t)-c'(t)\big\vert
  \convinfty{n} 0.
\end{align}
Then, in $\L^1(P),$
\begin{equation*}%\label{ShMcMCurve1}
  \frac{1}{\len\pi_n}\,\I(P_{L_{n}^{\pi}})\convinfty{n}
  \frac{1}{\len c}\int_0^T  h_{c'(t)}(P)\,dt.
\end{equation*}
\end{thm}

\pf
We have to show that
\begin{equation}\label{Curve1}
\lim_{n\to\infty}
\bigg\|\frac{1}{\len\pi_n}\I(P_{L_{n}^{\pi}})-
     \frac{1}{\len c}\int_0^T h_{c'(t)}(P)\,dt\,\bigg\|_{\L^1(P)}=0.
\end{equation}
Without loss of generality we can assume that $c$ has no self-intersections and that $c$ 
is parametrized by arc length.  As can be seen by the construction of the entropy for polygons,
\begin{equation*}
  \Big\|\frac{1}{\len\pi_n}\,\I(P_{L_{n}^{\pi}})-
  h_{\pi_n}(P) \Big\|_{\L^1(P)}
\end{equation*}
converges to $0.$
By the representation
formula in Remark \ref{IntReprhpi} and since
$\pi_n'(t)=((B_n\pi_n)^{-1})'(t/n),$
for all $t\in [0,nT],$ we obtain
\begin{align*}
h_{\pi_n}(P)= \frac{1}{\len\pi_n}\int_0^{nT}h_{\pi_n'(r)}(P)\,dr
            = \frac{n}{\len\pi_n}\int_0^T h_{((B_n\pi_n)^{-1})'(t)}(P)\,dt,
\end{align*}
and by (\ref{PolyCond}) and Lemma \ref{F_icontin}, the integral
converges to $ \int_0^T h_{c'(t)}\,dt.$ Use $\|\,\cdot\,\|$ for
the euclidian norm in the plane.  Using
\begin{equation*}
  \frac{1}{n}\len\pi_n =\int_0^T\|\pi_n'(nt)\|\,dt
  =\int_0^T\|(B_n^{-1})'(t)\|\,dt,
\end{equation*}
we obtain by (\ref{PolyCond})
\begin{equation*}
  \lim_{n\to\infty}\frac{1}{n}\len\pi_n
  =\int_0^T\|c'(t)\|\,dt
  =\len c.
\end{equation*}
This implies that
\begin{equation*}
  \Big\| h_{\pi_n}(P)-
  \frac{1}{\len c}\int_0^T h_{c'(t)}(P)\,dt \Big\|_{\L^1(P)}
\end{equation*}
converges to $0$ as well, and \eqref{Curve1} follows by the triangle
inequality.
\qed

Note that the limits do not depend on the sequence of polygons we
used to approximate the curve and that any approximation of the curve by lattice points 
can be described by a lattice approximation of a suitable polygon.  
This justifies the

\begin{definition}\label{DefEntrC}
Let $P$ and be as in Theorem \ref{ShMcMCurve}. Then
\begin{equation*}
  h_c(P):=\frac{1}{\len c}\int_0^T h_{c'(t)}(P)\,dt
\end{equation*}
is called {\rm specific entropy of $P$ along $c.$} 
\end{definition}

Note that the condition that the trace of $c$ does not contain the origin is no
real restriction for the definition as the expression on the left-hand side only
depends on the derivative of $c.$ This is intuitive, because we assumed that $P$
is stationary.  Note the following property for the entropies of the
blowups of a curve defined in \eqref{BlowingupDef}.  
The proof is a simple scaling argument.
\begin{cor}\label{Scaling}
Let $c:[0,T]\longrightarrow\R^2$ be a piecewise differentiable, and let 
$B_{\eta}c:[0,\eta T] \longrightarrow\R^2$ with 
$B_{\eta}c(t)=\eta\,c(\frac{t}{\eta})\ (\eta>0)$
be the family of its blowups. Then
%\begin{equation*}
$h_{B_{\eta} c}(P)= h_{c}(P)$ for all $\eta>0.$
%\end{equation*}
\end{cor}

%UUUUUUUUUUUUUUUUUUUUUUUUUUUUUUUUUUUUUUUUUUUUUUUUUUUUUUUUUUUUUUU
\subsection{Contour approximation}\label{sec:contours}
%UUUUUUUUUUUUUUUUUUUUUUUUUUUUUUUUUUUUUUUUUUUUUUUUUUUUUUUUUUUUUUU

In statistical mechanics, a {\it contour} is a set of sites
corresponding unambiguously to a chain of bonds.
Note that the lattice approximation is not a contour in this sense.
The last site before a new step is catercornered from the first site 
of the step, so not connected by a bond, and closing the gap
is not uniquely defined. We define the {\it contour approximation}
by adding, at each new step, the site which is one unit below it.
A new step begins in $i+1$ if and only if
$  \t_{\lambda}^i(\{a\})\ge 1-\lambda,$
and the site we add in this case is
$L_{\lambda,a}(I)-(0,1).$ 
For $z\in\Z$ and $u\in\N,$
$$  
\Lc_{\lambda,a}(z,\dots,z+u)
     :=L_{\lambda,a}(z,\dots,z+u) \cup
     \big\{L_{\lambda,a}(i)-(0,1)\,\big\vert\,
     0\le i\le u-1\;\wedge\;
    \t_{\lambda}^{z+i}(\{a\}) \ge 1-\lambda\big\} 
$$
be the contour approximation of the line segment $l_{\lambda,a}(I).$
Lemma \ref{laLa} and Lemma \ref{F_icontin} translate immediately to $\Lc.$
With little modifications we can prove contour versions of the 
limit theorems derived earlier in this section.  We will sketch the results here and
refer to \cite{Bre01SMcM} for details and proofs.  

For $n\in\N,$ define
$ \Lc_n(a):= \Lc_{\lambda,a_n}(z_n,\dots,z_n+u_n),$
and for $a\in\R$ and $z,{\widetilde z}\in\Z$ with ${\widetilde z}\ge z$ define
\begin{align}\label{Lf}
  \Lf_{\lambda,a}&({\widetilde z}, \dots,z):=\\ \nonumber
&\begin{cases}
   \Lc_{\lambda,a}({\widetilde z}-1,\dots,z)\cup
   \big(L_{\lambda,a}({\widetilde z})-(0,1)\big)
   &\qquad \mbox{{\rm if }}\t_{\lambda}^{{\widetilde z}-1}(\{a\})\ge 1-\lambda,\\
   \Lc_{\lambda,a}({\widetilde z}-1,\dots,z) &\qquad {\rm otherwise.}
\end{cases}
\end{align}
The {\rm specific contour entropy along a line}   $\widehat{h}_{\lambda}(P)$ 
with slope $\lambda$ is defined as
\begin{align*}
  \frac{1}{1+\lambda}\bigg(\int_0^1
  E\big[
    H\big( P_0[\cdot\,\vert\,
                           \F\big(\Lc_{\lambda,t}^{\sharp}(-\N)\big)
                       ](\w)\big)\big]\, dt
  +\int_{1-\lambda}^1
  E\big[
    H\big( P_0[\cdot\,\vert\,
			\F\big(\Lc_{\lambda,t}(-\N)\cup\{(0,1)\}\big)
                   ](\w)\big)\big]\, dt\bigg).
\end{align*}
We can show a contour version of Theorem \ref{ShMcMStretchL}:
\begin{equation*}
   \frac{1}{\vert \Lc_n(a)\vert}\,\I(P_{\Lc_n(a)})\convinfty{n} \widehat{h}_{\lambda}(P)
  \qquad\mbox{in }\L^1(P) \mbox{  and uniformly in } a\in\R.
\end{equation*}
The formulations of contour versions of Theorem \ref{ShMcMpolygon} 
and Theorem \ref{ShMcMCurve} are now obvious.

%UUUUUUUUUUUUUUUUUUUUUUUUUUUUUUUUUUUUUUUUUUUUUUUU
\section{Gibbs measures and specific entropies}\label{gib}
%\markboth{Large deviations}{Gibbs measures}
%UUUUUUUUUUUUUUUUUUUUUUUUUUUUUUUUUUUUUUUUUUUUUUUU

A collection $(U_V)_{V\subset\Z^d \text{ finite}}$ of functions on $\W$ 
is called {\it stationary summable interaction potential} if the following three
conditions are fulfilled:
(i) $U_V$ is measurable with respect to $\F_V$ for all $V\subset\Z^d.$
(ii) For all $i\in\N$ and all finite $V\subset\Z^d,$ $U_{V+i}=U_V\circ\theta_i.$
(iii) $\sum_{V\subset\Z^d\text{ finite}:\,0\in V} \parallel U_V\parallel_{\infty}<\infty.$
Let $\xi,\eta\in\W$ be two configurations. 
The {\it conditional energy} of $\xi$ on $V$
given the environment $\eta$ on $\Z^d\setminus V$ is defined as
\begin{equation*}%\label{ConEner}
  E_V(\xi\vert\eta)=
  \sum_{A\subset\Z^d\text{ finite}:\,A\cap V\neq\emptyset}
  U_A\big((\xi,\eta)_V\big),
\end{equation*}
where $(\xi,\eta)_V$ is the element of $\W$ given by
$(\xi,\eta)_V(i):=\xi(i),$ for $i\in V,$ and $(\xi,\eta)_V(i):=\eta(i)$ for 
$i\in \Z^d\setminus V.$
$P$ is called {\it Gibbs measure} with respect to $U$
if for any finite subset $V$ of $\Z^d$ the conditional distribution
of $\w_V$ under $P$ with respect to $\F_{\Z^d\setminus V}$ is given by
\begin{equation*}%\label{GibbsConDi}
P[\w_V=\xi_V\,\vert\,\F_{\Z^d\setminus V}](\eta)
=\frac{1}{Z_V(\eta)}e^{-E_V(\xi\vert\eta)},
\mbox{ where }
  Z_V(\eta):=\int_{\W}e^{-E_V(\xi\vert\eta)}\,P(d\xi)
\end{equation*}
is called {\it partition function}.
We say that there is a {\it phase transition} if there is more than one Gibbs
measure with respect to the same interaction potential.

Assume that $\Upsilon$ is furnished with a total order $\le,$ and
denote by $-$ the minimal and by $+$ the maximal element in $\Upsilon.$
Suppose that $U$ is {\it attractive} with respect to the order on $\Upsilon,$
in the sense of (9.7) in \cite{Pre76}.
Let $P^-$ and $P^+$ denote the minimal and the maximal Gibbs
measure with respect to $U,$ and let
$P^{\alpha}=\alpha P^-+(1-\alpha) P^+$  $(0<\alpha< 1)$
be their mixtures.
Both $P^-$ and $P^+$ are ergodic and, as follows from
\cite{Foe80}, they fulfill the strong $0$-$1$ law and the global Markov property.
F{\"o}llmer and Ort \cite{FoOrt88} define the {\it specific relative entropy}
based on hyperspaces by  
\begin{equation}\label{SpecSEntrFO}
s(P^-,P^+)=\frac{1}{d}\sum_{l=1}^d \int_{\W} H\Big(
        P_0^-\big[\cdot\big\vert\F^{(l)}\big](\w),
        P_0^+\big[\cdot\big\vert\F^{(l)}\big](\w)\Big)\,P^-(d\w),
\end{equation}
where $\F^{(l)}$ is the $\sigma$-algebra generated by those
coordinates in $\{(i^{(1)},...,i^{(d)})\in\Z^d\,\vert\,i^{(l)}=0\}$
which precede $0$ in the lexicographical order on $\Z^d.$

In the two-dimensional case, the conditions in \eqref{SpecSEntrFO} 
are simply along the coordinate axes.
Based on the work in Section \ref{gs}, we can now extend this definition to a
surface-order entropy along any direction $v\in S^1.$ Furthermore, we can
introduce an entropy along curves.  

\begin{definition}\label{RelEntrDir}
Let $v\in S^1.$  Let $c:[0,T]\longmapsto \R^2$ be a piecewise differentiable curve
parametrized by arc length with right derivative $c'.$
$$
  h_v(P^{-},P^{+}):=\int_0^1\int_{\W}
     H\Big( P_{0}^{-}\big[\cdot\big\vert\P_{\lambda(v),t}\big](\w),
           P_0^+\big[ \cdot\big\vert\P_{\lambda(v),t}\big](\w)
     \Big)\,P^{-}(d\w) dt
$$
is called {\rm specific relative entropy} of $P^{-}$ with respect to $P^{+}$ {\it in direction v.}
$$
h_c(P^{-},P^{+}):=\frac{1}{\len c}
\int_0^T h_{c'(t)}(P^{-},P^{+})\,dt
$$
is called {\rm specific relative entropy} of $P^{-}$ with respect to $P^{+}$ {\rm along} $c.$
\end{definition}

The order on $\Upsilon$ induces an order on the set $\M_1(\Upsilon)$
of probability measures on $\Upsilon:$ We say that $\mu$ is {\it larger} then
${\nu}$ if the density $\frac{d\mu}{d\nu}$ is an increasing
function with respect to the order on $\Upsilon,$ and in this case
we write $\mu\ge\nu.$
In particular, $\nu$ is absolutely continuous with respect to $\mu.$
The following inverse triangle inequality for relative entropies was shown
in the proof of Theorem 4.2 in \cite{FoOrt88}.  
For the reader's convenience
we state it in the following form.

\begin{lemma}\label{InvTriangle}
Let $\lambda\ge\mu\ge\nu,$ and assume that $\mu$ is bounded below
by a positive constant. Then
%\begin{equation*}
$H(\nu,\lambda)\ge H(\nu,\mu) + H(\mu,\lambda).$
%\end{equation*}
\end{lemma}

%UUUUUUUUUUUUUUUUUUUUUUUUUUUUUUUUUUUUUUUUUUUUUUUUUUUUUUUUUUU
\section{Lower bound}\label{lob}
%\markboth{Large deviations}{Lower bound}
%UUUUUUUUUUUUUUUUUUUUUUUUUUUUUUUUUUUUUUUUUUUUUUUUUUUUUUUUUUU

Let $P^{-}$ and $P^{+}$ be the minimal and maximal Gibbs measure defined in Section \ref{gib}.
The following lower bound for the large deviations of the {\it empirical field} 
$R_{n}(\w):=\sum_{i\in V_{n}}\delta_{\theta_{i}\w}$ of $P^+$ was proved by 
F{\"o}llmer and Ort \cite{FoOrt88}.
Recall that  $V_n$ is the set of all lattice sites in $[-n,n]^d,$ 
and that the boundary of a subset $V$ of $\Z^d$ is defined as
%\begin{align}\label{BdWh}
 $\partial V =\big\{i\in\Z^d\setminus V\,\big\vert\, \dist(i,V)=1 \big\}.$
%\end{align}

\begin{thm}\label{LowerBdFO}
For any open $A\in\M_1(\W),$
\begin{equation*}
  \liminf_{n\to\infty}\frac{1}{\vert\partial V_n\vert}\log
  P^{+}\big[R_n\in A\big]\ge -\inf_{\alpha:P_\alpha\in A}
  \sqrt{\alpha} \,s(P^-,P^+).
\end{equation*}
\end{thm}

The aim of this section is to improve the lower bound by
replacing the boxes by more general shapes in the two-dimensional case.  
The corresponding Shannon-McMillan theorems developed in Section \ref{gs}
will be the key to the proof.  
For a closed curve $c$ let $\inter c$ be the subset of
$\R^2$ surrounded by $c.$ Define the sets
\begin{align}\label{C_alpha}
  C_{\alpha}:=\big\{c\,\big\vert\,&
  c:[0,T]\longrightarrow\R^2\
    \text{closed piecewise
    $C^{1}$-curve parametrized by arc}\\
\nonumber
  &\text{length, without self-intersections, and with}\
  \area\inter c=\alpha\big\}.
\end{align}

%Recall that for a polygon $\pi$ $\pi_r, \lambda_r...$
%***********************************************************
%********************** Lower bound theorem ****************
\begin{thm}\label{LowerBd}
%Assume that $P^+$ has the local Markov property. Then
For any open $A\in\M_1(\W),$
\begin{equation*}
  \liminf_{n\to\infty}\frac{1}{\vert\,\partial V_n\,\vert}
  \log P^{+}\big[R_n\in A\big]
  \ge -\inf_{\alpha:P_\alpha\in A}
  \inf_{ c\in\C_{\alpha} }\frac{1}{4}
  \int_{0}^{T}\frac{dt}{\sqrt{1+\lambda(c'(t))^{2}}}
  \,\,h_c(P^-,P^+).
\end{equation*}
\end{thm}

\begin{rem}\label{SquareCurves}\rm{
Replacing the class $C_{\alpha}$ by squares with area $\alpha$
this bound coincides with the bound in Theorem \ref{LowerBdFO}:
Let $\pi$ be a square parametrized by arc length and with
$\area\inter \pi=\alpha.$ Then the length of every edge is
$\sqrt{\alpha}.$ For the two horizontal edges of the square the
slope $\lambda$ (cf.~\eqref{lambda_v}) is $0$ with respect to
the $x$-axis, and for the vertical edges it is $0$ with respect
to the $y$-axis. Therefore, the integral equals $4\sqrt{\alpha.}$
The entropy $h_{\pi}(P^-,P^+)$ equals $s(P^-,P^+),$ since the
$\sigma$-algebras $\P_{0,t}$ coincide with $\F^{(2)}$ for the
horizontal edges and with $\F^{(1)}$ for the vertical edges.}
\end{rem}

\begin{rem}\label{LowerBc}\rm{
In the case where the Markov property holds only with respect to the contour boundary 
we can state a bound similar to the one in Theorem \ref{LowerBd} by replacing the 
lattice approximation by the contour aproximation.
The proof is essentially the same as for Theorem \ref{LowerBd}
(cf.~\cite{Bre01SMcM} for details).}
\end{rem}

%--------------------- Preparation of proof (geometric lemmata) ------------------------
The rest of this section is devoted to the proof of Theorem \ref{LowerBd}.  
We will need some properties of the geometry of lattice approximations of polygons 
and their interplay with the random field.
To begin with, we restate explicitly the global Markov property  for
random fields in the case when the conditioning is concentrated on a
set of sites surrounded by a closed polygon $\pi$ without self-intersections. 
We use the notation
%\begin{equation}\label{Gammac}
$\Gamma(c):=\Inter c\cap\Z^2$
%\end{equation}
to indicate the set of lattice points surrounded by a closed curve $c.$  
By the definition \eqref{LBnpi}, $\partial \big(\Z^2\setminus\Gamma(\pi)\big)=L^{\pi},$ 
and the global Markov property (cf.~Definition \eqref{MP}), with
$V=\Z^2\setminus\Gamma(\pi),$ we have that for any
$\F(\Z^2\setminus\Gamma(\pi))$-measurable nonnegative function $\Phi,$
\begin{equation}\label{MPIntPolyg}  
E\big[\phi\,\big\vert\,\F_{\Gamma(\pi)}\big]=
E\big[\phi\,\big\vert\,\F_{L^{\pi}}\big].
\end{equation}
We will further need two lemmata that compute
the asymptotic fractions of the lengths of the blowups of a line segment,
or a polygon, and the sizes of their lattice approximations.

\begin{lemma}\label{LovB}
Let $I$ be a real interval,
$l(x)=\lambda x+a$ be a linear function with slope  $\lambda,$ and
$B_{k}\ (k\in\N)$ be the sequence of its blowups restricted to $I.$
If $L_{k}$ is the lattice approximation of $B_{k}$ then
\begin{equation*}
\lim_{k\to\infty}\frac{\vert\,L_{k}\,\vert}{\len B_{k}}
=\frac{1}{\sqrt{1+\lambda^{2}}}.
\end{equation*}
\end{lemma}
\pf
We consider only the case when $0\le\lambda\le 1,$
that is, when the lattice approximation is given by $L(z)=(z,[l(z)])$
$(z\in I\cap\Z).$
Other cases only differ in terms of notation.
For any $k\in\N,$ $\vert L_{k}\vert$ is either $[\len b_{k}]$ or
$[\len b_{k}]+1,$ where $b_{k}$ is the projection of $B_{k}$ to
the $x$-axis. We can ignore the second case, since the additional
point does not matter for the limit. Observe
that $ (\len B_{k})^{2}=(\len b_{k})^{2}+(\lambda\, \len
b_{k})^{2}. $ Consequently,
%\begin{equation*}
$\len b_{k}=\len B_{k}\,/\,\big(\sqrt{1+\lambda^{2}}\big),$
%\end{equation*}
which proves the convergence.
\qed

\begin{lemma}\label{LovBPi}
Let $\pi$ be a polygon with edges $\pi_r$ $(r=1,\dots,R)$
and $\lambda_r$  $(r=1,\dots,R)$ their slopes as
defined in \eqref{lambda_v}.
Let $B_k\pi\ (k\in\N)$ be the blowups of $\pi$ and
$L_k\pi\ (k\in\N)$ their lattice approximations.
Then
\begin{equation*}%\label{LovBLimPi}
   \lim_{k\to\infty}\frac{\vert\,L_{k}\pi\,\vert}{\len B_{k}\pi}
   =\sum_{r=1}^R\frac{1}{\sqrt{1+\lambda_r^{2}}}\frac{\len\pi_r}{\len\pi}.
\end{equation*}
\end{lemma}
\pf
Using $\len B_k\pi_r=k\,\len\pi_r,$ we obtain
$$
\frac{\vert\,L_{k}\pi\,\vert}{\len B_{k}\pi}
=\sum_{r=1}^R\frac{\vert\,L_{k}\pi\,\vert}{\len B_{k}\pi_r}
\frac{\len B_{k}\pi_r}{\len B_{k}\pi}
=\sum_{r=1}^R\frac{\vert\,L_{k}\pi\,\vert}{\len B_{k}\pi_r}
\frac{\len\pi_r}{\len \pi}.
$$
By the previous lemma applied to the individual sides,
the first factors converge to $1\big/\sqrt{1+\lambda_r^{2}},$
which prove the statement of the lemma.
\qed

%################### PROOF OF THEOREM #################
%E\big[\phi\,\vert\,\F_{\partial(\Z^2\setminus\Gamma(\pi))}\big]=
{\bf Proof of Theorem \ref{LowerBd}:}
Let be $0<\alpha\le 1,$ such that $P_{\alpha}\in A.$ Since $A$ is open, we can
choose open neigborhoods $A^-$ and $A^+$ of $P^-$ respectively
$P^+$ in $\M_1(\W)$ such that
$
\alpha A^- + (1-\alpha) A^+\subseteq A.
$
Without loss of generality we may assume that $A^-$ and
$A^+$ are in $\F_{V_p}$ for some $p\in\N.$
%$\MMM$\\
Define the set
%\begin{equation*}%\label{Pi_alpha}
$
\Pi_{\alpha}:=\big\{\pi\,\big\vert\,\pi\ \text{closed polygon without
self-intersections, }
\area\inter\pi=\alpha\big\}.
$
%\end{equation*}
Let $\pi\in\Pi_{4\alpha}$ with $0\in\inter\pi,$ and let
$(B_n\pi)_{n\in\N}$ be the sequence of blowups of $\pi.$
For $\alpha=1$ take $C_n:=V_n.$ Otherwise, define
\begin{equation*}%\label{CnDn}
    C_n:=\Gamma(B_{k(n)}\pi)\
    \quad\mbox{and}\quad
    D_n:=V_n\setminus\Gamma(B_{l(n)}\pi),
\end{equation*}
where $k(n)$ and $l(n)$ are chosen such that $k(n)\le l(n),$
%\begin{equation*}
$l(n)-k(n)\convinfty{n}\infty,$
%\end{equation*}
\begin{equation}\label{wantkl}
  \lim_{n\to\infty}\frac{\vert\, C_n\,\vert}{\vert\, V_n\,\vert}=\alpha
  \qquad\text{and}\qquad
  \lim_{n\to\infty}\frac{\vert\, D_n\,\vert}{\vert
  V_n\,\vert}=1-\alpha.
\end{equation}
To see that such sequences exist we show that 
$k(n)=\big[\sqrt{ {\alpha\vert V_n\vert}/{\Area\Inter \pi}\, }\,\big]$
and
$l(n)=\big[k(n)+\sqrt{n}\,\big]$ fulfill the conditions. 
Obviously, both the sequences and their difference tend to infinity
as $n$ goes to infinity. Using 
$\Area\Inter(B_k\pi)=k^2\Area\Inter \pi,$
\begin{equation}\label{wantkl1}
  \frac{\vert\,\Gamma(B_k\pi)\,\vert}{\Area\Inter B_k\pi}
  \convinfty{k} 1\qquad\mbox{and}\qquad
\frac{k(n)^2}{\vert\,V_n\,\vert}\convinfty{n}\frac{\alpha}{\Area\Inter\pi}.
\end{equation}
We obtain for the first expression in (\ref{wantkl})
$$
\lim_{n\to\infty}\frac{\vert\, C_n\,\vert}{\vert\, V_n\,\vert}
=\lim_{n\to\infty}\frac{\Area \Inter(B_{k(n)}\pi)}{\vert\,V_n\,\vert}
=\lim_{n\to\infty}\frac{k(n)^2\,\Area\Inter\pi}{\vert\,V_n\,\vert}
=\alpha.
$$
Similarly we see for the second expression in \eqref{wantkl},
$$
\lim_{n\to\infty}\frac{\vert\, D_n\,\vert}{\vert\, V_n\,\vert}
=1-\lim_{n\to\infty}\frac{\vert\,\Gamma(B_{l(n)})\,\vert}{\vert\,V_n\,\vert}
= 1- \lim_{n\to\infty}\frac{l(n)^2\,\Area\Inter\pi}{\vert\,V_n\,\vert}.
$$
By definition of $l(n),$ $l(n)^2=k(n)^2+[2k(n)\sqrt{n}]+n.$
The last two addends are of order $n$ and will go to $0$ when
divided by $\vert\,V_n\,\vert.$   It remains to study 
$k(n)^2\,\Area\Inter\pi\,/\,\vert\,V_n\,\vert,$  but we already know from the 
second statement in \eqref{wantkl1} that this converges to $\alpha.$

Define
$$
R_n^-=\frac{1}{\vert\, C_{n,p}\,\vert}
\sum_{i\in C_{n,p}}\delta_{\theta_i}\w
\qquad\text{and}\qquad
R_n^+=\frac{1}{\vert\, D_{n,p}\,\vert}
\sum_{i\in D_{n,p}}\delta_{\theta_i}\w,
$$
where $C_{n,p}:=V_{k(n)-p}$ and
$D_{n,p}:=V_{n-p}\setminus V_{l(n)+p}.$ Then
%\begin{equation}\label{R_NinFC}
$\{R_n^-\in A^-\}\in\F_{C_n},$ $\{R_n^+\in A^+\}\in\F_{D_n},$
%\end{equation}
for large enough $n,$
%\begin{equation*}%\label{R_ninLambda_n}
  $\{R_n\in A\}\supseteq \{R_n^-\in A^-\}\cap \{R_n^+\in A^+\} := \Lambda_n.$
%\end{equation*}
Define the measures 
$$
Q_n= P_{C_n}^-\otimes P_{\Z^2\setminus
C_n}^{+}
\qquad(n\in\N).
$$ 
$Q_n$ coincides with $P^-$ on
$\F_{C_n}$ and with $P^+$ on $\F_{D_n},$ and makes these
$\sigma$-fields independent. Thus %, and by \eqref{R_NinFC}, 
we obtain $ Q_n[ \Lambda_n]=P^-[R_n^-\in A^-]\,P^+[R_n^+\in A^+], $
and by the ergodic behaviour of $P^-$ and $P^+,$ the sequence
%\begin{equation}\label{Q_nto1}
$Q_n[ \Lambda_n]$ $(n\in\N)$ converges to $1.$
%\end{equation}
Let $\phi_n$ denote the density of $Q_n$ with respect to $P^+$ on
$\F_{C_n\cup D_n}.$ Then for $\gamma>0, \varepsilon>0,$ and for
large enough $n,$
\begin{align*}
  P^+[R_n\in A]  \ge P^+[ \Lambda_n]
  & \ge \int 1_{ \Lambda_n\cap\{\vert\partial V_n\vert^{-1}
             \log\phi_n \le \gamma+\varepsilon \} }
  \phi_n^{-1}\,dQ_n\\
  & \ge \exp(-(\gamma+\varepsilon)\,\vert\partial V_n\vert)\,
  Q_n\big[ \Lambda_n\cap\big\{ \vert\partial V_n\vert^{-1}
             \log\phi_n \le \gamma+\varepsilon\big\}\big].
\end{align*}
By the convergence of the $Q_n[ \Lambda_n],$ the lower bound
%\begin{equation*}
  $\liminf_{n\to\infty}\vert\partial V_n\vert^{-1}
  \log P^+[R_n\in A]\ge -\gamma$
%\end{equation*}
follows if $\gamma$ is chosen such that, for any $\varepsilon>0,$
\begin{equation}\label{Q_nPhi_n}
  \lim_{n\to\infty}Q_n\big[ \vert\partial V_n\vert^{-1}
  \log\phi_n\le \gamma+\varepsilon\big]=1.
\end{equation}

We will show that (\ref{Q_nPhi_n}) holds with
$$
\gamma=
\sum_{r=1}^R\frac{1}{\sqrt{1+\lambda_r^{2}}}\frac{\len\pi_r}{8}\,
\,h_{\pi}(P^-,P^+).
$$
Since $Q_n=P^+$ on $D_n,$ and the fact that both $P^-$ and $P^+$
are Gibbs measures with respect to the same potential we obtain
\begin{equation*}
  \phi_n(\w)=\frac{P^-[ \w_{C_n}]P^+[ \w_{D_n}]}
  {P^+[ \w_{C_n\cup D_n}]}
  =\frac{P^-[ \w_{C_n}]}
  {P^-[\w_{C_n}\,\vert\,\w_{D_n}]}.
\end{equation*}
Let $L_n$ be the lattice approximation %(cf.~Definition \ref{DefL})
of $B_n\pi.$ By (\ref{MPIntPolyg}),
\begin{align*}
  &P^-[ \w_{C_n}\,\vert\,\w_{D_n}]
  =P^-[ \w_{D_n}\,\vert\,\w_{C_n}]\frac{P^-[ \w_{D_n}]}{P^-[ \w_{C_n}]}\\
  &=P^-[ \w_{D_n}\,\vert\,\w_{L_{k(n)}}]\frac{P^-[ \w_{D_n}]}{P^-[
\w_{C_n}]}
  =P^-[ \w_{L_{k(n)}}\,\vert\,\w_{D_n}]\frac{P^-[ \w_{L_{k(n)}}]}{P^-[
\w_{C_n}]},
\end{align*}
and thus
\begin{equation}\label{Phi}
  \phi_n(\w)=\frac{P^-(\w_{L_{k(n)}})}
  {P^-(\w_{L_{k(n)}}\,\vert\,\w_{D_n})}
  =\frac{P^-(\w_{L_{k(n)}})}
  {P^+(\w_{L_{k(n)}}\,\vert\,\w_{D_n})}.
\end{equation}
Going around the $R$ sides of $L_{k(n)}$, and conditioning
site by site as in the proof of Theorem \ref{ShMcMpolygon}, we obtain
$$
\frac{1}{\vert\,V_n\,\vert}\log\phi_n(\w)=
\frac{1}{\vert\,V_n\,\vert}\sum_{r=1}^R\Psi^{(r)},
$$
where the $\Psi^{(r)}$ corresponds to the $r$-th side of the
polygon. Similar to the calculation between \eqref{pfpoly1}
to \eqref{pfpoly4} we obtain
\begin{align*}
  \Psi^{(r)}=\sum_{i=0}^{u_n} & Z_{n,i,t}
  \circ\theta_{L_{\lambda,\t_\lambda^{z_n}(\{t_n\})}(i)}
           \circ\theta_{L_{\lambda,t_n}(z_n)},
\end{align*}
where $\lambda$ is the slope of the $r$th side of the polygon,
$t_{n}$ and $r_{n}$ are as in Subsection \ref{sec:poly}, and
%\begin{align*}
 $ Z_{n,i,t}= X_{n,i,t} -Y_{n,i,t},$
%\end{align*}
with
\begin{align*}
X_{n,i,t}=& \log P_0^{-}\Big(\w(0,0)\,\Big\vert\,
        \w\big(\Lf_{\lambda,\t_{\lambda}^{z_n+i}(\{t_n\}) }
        (-1,\dots,-i)\cup A_{n,i,t}\big)\Big),\\
\text{and}\ Y_{n,i,t}=& \log P_0^{+}\Big(\w(0,0)\,\Big\vert\,
        \w\big(\Lf_{\lambda,\t_{\lambda}^{z_n+i}(\{t_n\}) }
        (-1,\dots,-i)\cup B_{n,i,t}\big)\Big).
\end{align*}
To simplify notation we have omitted the index $r.$ For the sets in the
conditional expectations we have
$ A_{n,i,t}\subseteq
  B_{n,i,t}\subseteq
  \Z^2\setminus\big(H_n-L_{\lambda,t_n}(z_n+i)\big).$
$A_{n,i,t}$ is obtained by
shifting a subset of $L_n\subseteq C_n.$
$H_n$ is constructed as in the paragraph above, %(\ref{OutHex1})
%$MMM$
but using the minimum of the diameter $d_n$ and the distance
$l(n)-k(n)$ in place of $d_n.$

To prove convergence, we study the $X$ and $Y$-parts separately.
Because of the way the sets $A_{n,i,t}$ are constructed, the behavior of
$X_{n,i,t}$ under $Q_n$ is the same as under $P^-.$ But the proof
of Theorem \ref{ShMcMpolygon} shows that
\begin{align*}
  \sum_{i=0}^{u_n} & X_{n,i,t}
  \circ\theta_{L_{\lambda,\t_\lambda^{z_n}(\{t_n\})}(i)}
           \circ\theta_{L_{\lambda,t_n}(z_n)}
\end{align*}
converges to $-h_{\pi}(P^-)$ in $\L^1(P^-).$
The convergence remains true when we replace $X_{n,i,t}$ by
\begin{equation*}
        X_{n,i,t}^-:=\log P_0^{-}\Big(\w(0,0)\,\Big\vert\,
        \w\big(L_{\lambda,\t_{\lambda}^{z_n+i}(\{t_n\}) }
        (-1,\dots,-i)\big)^-\Big),
\end{equation*}
where, for a subset $L$ of $\Z^2,$ the element
%\begin{equation}\label{w_}
$\w\big(L\big)^{-}$
%\end{equation}
equals $\w$ on $L$ and assumes the minimal state in $\Upsilon$ outside of
$H_n-L_{\lambda,t_n}(z_n+i).$
To control the behavior of $Y_{n,i,t}$ under $Q_n$ define
$Z_{n,i,t}^-=X_{n,i,t}^{-} -Y_{n,i,t}.$
Use the law of large numbers for martingales with bounded
increments in its $\L^2$-form to replace
\begin{align*}
  \frac{1}{\vert\,L_{k(n)}\,\vert}
  \sum_{i=0}^{u_n} Z_{n,i,t}^{-}
  \circ\theta_{L_{\lambda,\t_\lambda^{z_n}(\{t_n\})}(i)}
           \circ\theta_{L_{\lambda,t_n}(z_n)}
\end{align*}
by
\begin{align*}
  \frac{1}{\vert\,L_{k(n)}\,\vert}
  \sum_{i=0}^{u_n} E\big[ Z_{n,i,t}^{-}
  \circ\theta_{L_{\lambda,\t_\lambda^{z_n}(\{t_n\})}(i)}
           \circ\theta_{L_{\lambda,t_n}(z_n)}\,\big\vert\,\A_{n,i,t}\big],
\end{align*}
where $\A_{n,i,t}$ is the $\sigma$-field generated by
the sites in $D_n$ and those sites of $L_{k(n)}$ which precede $i$ in
the canonical ordering of $L_{k(n)}.$ $\MMM$
These conditional expectations can be written as the relative
entropy $H(\nu,\mu),$ with the random measures
\begin{align*}
  \mu(\w):=&
  P_0^-\big[\,\cdot\,\big\vert\,
    \w\big(L_{\lambda,t_{\lambda}^{z_n+i}(\{t_n\})}
        (-1,\dots,-i)\big)^-\big]\\
\text{and}\qquad  \nu(\w):= &
  P_0^{+}\big[\,\cdot\,\,\big\vert\,
     \w\big(L_{\lambda,\t_{\lambda}^{z_n+i}(\{t_n\}) }
        (-1,\dots,-i)\cup B_{n,i,t}\big)\big].
\end{align*}

Now we want to replace $\mu$ by a measure $\eta$ for which
\begin{equation}\label{LamConv}
\big({\vert\,L_{\lambda,\t_{\lambda}^{z_n+i}(\{t_n\})}
        (-1,\dots,-i)\,\vert}\big)^{-1}
  \sum_{i=0}^{u_n} H(\nu,\eta)
  \circ\theta_{L_{\lambda,\t_\lambda^{z_n}(\{t_n\})}(i)}
           \circ\theta_{L_{\lambda,t_n}(z_n)}
\end{equation}
converges to $h_{\lambda}(P^{-},P^{+}),$
in $\L^1(P^{-}),$ as $n$ goes to infinity.
Define $\w(L)^{+}$ in analogy to $\w(L)^{-}.$ Since for {\it all} $\w,$
$$
P_0^-\Big[\,\cdot\,\Big\vert\,\w\big(L_{\lambda,t}
        (-1,\dots,-i)\big)^{+}\Big]\convinfty{i}
P_0^-\Big[\,\cdot\,\Big\vert\,\P_{\lambda,t}\Big],
$$
we obtain \eqref{LamConv} by taking
$$
\eta(\w):=P_0^-\Big[\,\cdot\,\Big\vert\,\w
\big(L_{\lambda,\t_{\lambda}^{z_n+i}(\{t_n\})}
        (-1,\dots,-i)\big)^{+}\Big].
$$
By Lemma \ref{InvTriangle},
$H(\nu(\w),\mu(\w))\le H(\nu(\w),\eta(\w)).$
Summing over $r=1,\dots,R,$ and
passing from convergence in $\L^1(P^-)$ to stochastic convergence
with respect to $Q_n$ yields
\begin{equation*}%\label{Q_nPhi_n2}
  \lim_{n\to\infty} Q_n\big[\vert\,L_{k(n)}\,\vert^{-1}
  \phi_n> h_{\pi}(P^{-},P^{+})+\varepsilon\big]=0
\end{equation*}
for any $\varepsilon>0.$
To derive \eqref{Q_nPhi_n} 
%with
%$$\gamma=
%\sum_{r=1}^R\frac{1}{\sqrt{1+\lambda_r^{2}}}\frac{\len\pi_r}{8}\,\,h_{\pi}(P^-,P^+),$$
it remains to show that
\begin{equation} \label{L_kovPV}
\lim_{n\to\infty}
\frac{\vert\,L_{k(n)}\,\vert}{\vert\,\partial V_n\,\vert}
=\sum_{r=1}^R\frac{1}{\sqrt{1+\lambda_r^{2}}}
\frac{\len\pi_r}{8}.
\end{equation}
%this is equivalent to \eqref{Q_nPhi_n} with $c=...$
The fraction on the left-hand side can be written as a product:
\begin{equation*}\label{ThreeF}
\frac{\vert\,L_{k(n)}\,\vert}{\vert\,\partial V_n\,\vert}=
\,\frac{\vert\,L_{k(n)}\,\vert}{\len B_{k(n)}\pi}
\cdot\frac{k(n)\,\len \pi}{\vert\,\partial V_n\,\vert}.
\end{equation*}
We first study the asymptotics of the second factor:
As we are only interested in the limit behaviour, 
we can drop the brackets for the integer part in the definition of $k(n)$ and just use
 $\sqrt{\alpha\vert\, V_n\,\vert/\,\area \inter \pi}.$
The denominator equals $4\,\vert\,V_n\,\vert.$  Reducing the fraction yields $1/8.$
The limit of the first factor in \eqref{ThreeF} was computed in Lemma \ref{LovBPi}.
So, the second factor in \eqref{ThreeF} converges to $\len \pi/8.$
Combining all yields \eqref{L_kovPV}.

Now, we replace the polygon $\pi$ by the polygon
$\widetilde{\pi}=B_{1/2}\pi.$ Since
%\area\widetilde{\pi}=\alpha\qquad\text{and}\qquad
$\len\widetilde{\pi}_r=\len\pi_r/2,$
and since, by Corollary \ref{Scaling},%
$h_{\widetilde{\pi}}(P^{-},P^{+})
=h_{\pi}(P^{-},P^{+}),$
%\end{document}
$$
\gamma
= \sum_{r=1}^R\frac{1}{\sqrt{1+\lambda_r^{2}}}
\frac{\len\widetilde{\pi}_r}{4}\,h_{\widetilde{\pi}}(P^-,P^+).
$$
Finally, by Lemma \ref{F_icontin}, the infimum of that function
over all polygons $\widetilde{\pi}\in\Pi_{\alpha}$
equals the infimum over all curves $c\in C_{\alpha}.$
%which yields the bound \eqref{LowerBd}.
\qed

\nocite{Mil88}
\nocite{DeS89}

%AAAAAAAAAAAAAAAAAAAAAAAAAAAAAAAAAAAAAAAAAAAAAAAAAAAAAAAA
\section*{Acknowledgements}
%AAAAAAAAAAAAAAAAAAAAAAAAAAAAAAAAAAAAAAAAAAAAAAAAAAAAAAAA

%\nd {\bf Acknowledgements.}
This work was a part of my PhD thesis, and I thank
my supervisor Hans F{\"o}llmer for suggesting to extend the concept of 
surface entropies to general shapes via a stepwise approach, 
and to use them to refine the large deviation lower bounds in the phase-transition regime. 
I am grateful to Hans-Otto Georgii for noticing and correcting a 
mistake in the constant in the lower bound, and I thank Dima Ioffe for interesting
discussions.  I also want to thank the referee for insightful and supportive advice.
Finally, I am very grateful to Kristan Aronson for asking
a very good question about this work.

%BBBBBBBBBBBBBBBBBBBBBBBBBBBBBBBBBBBBBBBBBBBBBBBBBBBBBBBBBBBBBBBB
%%\chapter{Bibliography}\nonumber
%%\chapter*{Bibliography}\nonumber
%\addcontentsline{toc}{chapter}{\numberline{}Bibliography}
%\bibliographystyle{plain}
%\bibliography{biblio}
%BBBBBBBBBBBBBBBBBBBBBBBBBBBBBBBBBBBBBBBBBBBBBBBBBBBBBBBBBBBBBBBB

\smallskip

\noindent
{\em
Julia Brettschneider\\[1.5mm]
Department of Statistics, University of Warwick, Coventry, CV4 7AL, UK\\[1.5mm]   
Department of Community Health $\&$ Epidemiology and Cancer Research Institute Division of \\
Cancer Care $\&$ Epidemiology, Queen's University, Ontario, K7L 3N6, Canada  \\[1.5mm]    
{\tt  julia.brettschneider@warwick.ac.uk}               
}

\end{document}